\documentclass [11pt] {article}
\usepackage {latexsym}
\usepackage {amssymb}
\usepackage[dvips]{graphicx}

\newcommand {\zz} {\mathbb Z}
\newcommand {\mycirc}{*}

\newtheorem{theorem}{Theorem}
\newtheorem{lemma}{Lemma}[section]
\newtheorem{definition}{Definition}[section]

\newtheorem{remark}{Remark}[section]

\newcounter{intr}
\renewcommand{\theintr}{{\bf \thesection.\arabic{intr}}}
\newenvironment{intr}{\noindent \theintr \refstepcounter{intr}}{\medskip}
\newenvironment{proof}{\noindent{\bf Proof.}}{$\Box$\par\medskip}
\newenvironment{vtheorem}[2]{\noindent{\bf #1 #2.}\em }{\par\medskip}

\begin {document}
\begin {center}
{\Large\bf  The topological classification of Morse-Smale diffeomorphisms on 2-manifolds.}

\vspace {3mm}
{\bf Igor Vlasenko}\\
Institute of Mathematics, \\
Tereshchenkivska, 3\\
Kiev, Ukraine. \\
e-mail: vlasenko@imath.kiev.ua
\end {center} 

\begin{center}
{\bf Abstract}
\end{center}
\begin{quotation}
\small
We consider general Morse-Smale diffeomorphisms on a closed orientable 
two-dimensional surface. 
We prove that topological invariant which uniquely determine Morse-Smale diffeomorphisms up to the topological conjugacy is finite. It is obtained in an explicit form that gives us the topological classification of Morse-Smale diffeomorphisms.
\end{quotation}


This paper is devoted to the problem of the topological classification of Morse-Smale diffeomorphisms on two-dimensional manifolds. It presents a complete finite topological invariant of general Morse-Smale diffeomorphisms and necessary and sufficient conditions of topological conjugacy of Morse-Smale diffeomorphisms. 

The basic difficulty here is presence of a complicated structure of intersection of stable and unstable manifolds of periodic points, generating in the general case an infinite number of periodic trajectories. This problem has a rather large bibliography (see \cite{Bezen} -- \cite{Bor5}). 


The most general result was obtained by A. Z. Grines, who solved  the problem for the case of a finite number of heteroclinic trajectories in the paper \cite{Gri1}. Different special cases were investigated in papers \cite{Bgri1}--\cite{Bgri4}. The topological classification of periodic components of diffeomorphisms was given by A. Bezen in paper \cite{Bezen}. 
This problem was also considered in papers \cite{Bor1}, \cite{Bor2}, but they require 
essential refinements.

This paper considerably improves the initial author's proof that first appeared  in two articles \cite{Vla1}, \cite{Vla2} (in Russian).

  The special gratitude the author expresses to my supervisor professor Vladimir V. Sharko and to my colleagues Serge Maksimenko, Mark Pankov, Eugene Polylyakh and Alexander Prishlak.

\section {Preliminary information.} \par
\begin{intr}
Let $f$ be a homeomorphism on a compact manifold $M$. 
A point $x \in M$ is a {\sl nonwandering point} of $f$ provided  $\forall U(x)$ $\exists m \not= 0$ that $f^{m}(U) \cap U \not= \varnothing$. (We denote by $U(x)$ a neighborhood of a point $x$). 
Otherwise point is called a \textit{wandering} one. The set of all nonwandering points $f$ is denoted by $\Omega (f)$. The point set $\{f^{k}(x),  k \in Z \}$ is called \textit {a trajectory of the point $x$}. 
A point $x$ is called \textit{periodic} with a period $n$ if $f^{n}(x)= x$ and $\forall k = 1... n-1$ $f^k(x) \not= x$. The set of all periodic points $f$ is denoted by $Per (f)$. \par
\end{intr}

\begin{intr}\label{intr12}
Let $Diff(M)$ be a space of diffeomorphisms on a two-dimensional compact manifold $M$, $f$ be a diffeomorphism. 
The point $x \in Per(f)$ of a period $m$ is called \textit{hyperbolic} if the derivative $D_{x}f^{m}$ of $f^{m}$ at $x$ has the spectrum disjoint from the unit circle. If all the eigenvalues of $D_{x}f^{m}$ lie inside the unit circle, $x$ is called \textit{a sink}. When each of the eigenvalues has modulus greater than one, $x$ is called  \textit{a source}. Otherwise $x$ is called \textit{a saddle}. 
We denote a set of saddle periodic points of  $f$ by $\Omega'(f)$. 
Let $\Gamma$ be a hyperbolic periodic trajectory with a period $m$ and $x \in \Gamma$. Let $E^{u}$ denotes a subspace $T_{x}M$, spanned by the eigenvalues of the $D_{x}f^{m}$, corresponding to eigenvalues whose moduli are greater than one. Let $E^{s}$ be a subspace $T_{x}M$, spanned by the remaining eigenvalues. Define the orientation type of the trajectory $\Gamma$ Will be $+1$, if $D_{x}f^{m}:E^{u}\rightarrow E^{u}$ preserves the orientation and $-1$ if it reverses the orientation. 
\end{intr}

\begin{intr}
Let $d$ be a metric on $M$, $x$ be a hyperbolic fixed (i. e. $f (x)=x$) point of $f$.
 The stable and unstable manifolds of  $x$ are: \par
 $$W^{s}(x) =\{y \in M\: |\: d (x, f (y)) \rightarrow 0, n \rightarrow + \infty \}$$
 $$W^{u}(x) =\{y \in M\: |\: d (x, f (y)) \rightarrow 0, n \rightarrow- \infty \}$$
  For a hyperbolic point $x$ of a period $m$ the stable and unstable manifolds are defined to be the stable and unstable manifolds of $x$ under $f^{m}$. \par

For a hyperbolic trajectory the stable and unstable manifolds are defined to be the union of the stable and unstable manifolds of points belonging to the trajectory.

We also introduce the following useful designations: 
We denote by $W_{i}^{\sigma}(x)$, where $\sigma \in \{u, s \}$, $i = 1,2$,  pathcomponents of $W^{\sigma}(x) \setminus \{x \}$. Here $\sigma$ denotes either $u$ or $s$ and $\rho$ denotes a manifold dual to $\sigma$ : if $ \sigma = u $ then $ \rho = s$ and vice versa. The set of saddle periodic points of $f$ is denoted by $Per'(f)$.  Also
$$
\partial W^{\sigma}(x)=\textbf{cl}W^{\sigma}(x) \setminus W^{\sigma}(x)\qquad  \partial W_{i}^{\sigma}(x) =\textbf{cl}W_{i }^{\sigma}(x) \setminus W_{i}^{\sigma}(x)
$$
 where $\textbf{cl}$ denotes closure, 
$$
W'{}^{u}=\bigcup_{x\in Per'(f)}W^{u}(x) \qquad W'{}^{s}=\bigcup_{x \in Per'(f)}W^{s}(x)\qquad W'= W'{}^{u}\bigcup W'{}^{s}
$$
The points of  intersection $W'{}^{u} \cap W'{}^{s} $ are called \textit{heteroclinic points}. 
\end{intr}

\begin{intr}
$f$ is called a Morse-Smale diffeomorphism, if 

1. $ \Omega (f) $ is finite

2. all periodic points are hyperbolic

3. for all $x, y\in \Omega (f)$ $W^{u}(x)$ and $W^{s}(y)$ have transversal intersection. 

\medskip\noindent
Further in the paper $f$ will denote everywhere a $C^r$ Morse-Smale diffeomorphism. \par
\end{intr}

\begin{intr}
Two diffeomorphisms $f$ and $g$ are topologically conjugate if and only if there exist a homeomorphism $h: M \rightarrow M$ such that $ h \circ f = g \circ h $. 

A diffeomorphism $f$ is structurally stable provided there exists a neighborhood $U$
of $f$ in $Diff(M)$ such that there exists a homeomorphism $h: M \rightarrow M$ such that $h \circ f = g \circ h $ for each $h \in U$.

J. Palis and S. Smale proved in \cite{ps} that Morse-Smale diffeomorphisms are exactly the structurally stable diffeomorphisms with a finite nonwandering set.
\end{intr}

\begin{intr}\label{intr13}
We list here some properties of Morse-Smale diffeomorphisms from Smale's papers  \cite{S1},\cite{S2}: 
\begin{enumerate}
\item $ \forall x \in \Omega (f)\quad \partial W^{u}(x) \subset W'{}^{u} \cup \Omega^{s}(f)$, $\partial W^{s}(x) \subset W'{}^{s} \cup \Omega^{u}(f)$, where $\Omega^{u}(f)$, $\Omega^{s}(f)$ are the unions of sinks and sources respectively. \par

\item If $\partial W^{u}(x_{i}) \cap W^{u}(x_{j}) \not= \varnothing$ then $W^{u}(x_{i}) \cap W^{s}(x_{j}) \not= \varnothing $. \par

\item If $W^{u}(x_{i}) \cap W^{s}(x_{j}) \not= \varnothing$ then $\partial W^{u}(x_{i}) \supset W^{u}(x_{j}) $. \par

\item If $W^{u}(x) \cap W^{s}(y) \not= \varnothing$ and $W^{u}(y)\cap W^{s}(z)\not=  \varnothing $\\ then  $W^{u}(x) \cap W^{s}(z) \not= \nobreak \varnothing $. \par

\item All $W^{\sigma}(x)$ $x \in \Omega (f) $ are immersed submanifolds $ M $. \par
\end{enumerate}
\end{intr}

\begin{intr}
A domain $S \subset M$ is called \textit{a domain of a wandering type}, if $f^{k} (S) \cap S = \varnothing$, $k \not= 0$. An area $S \subset M$ is called \textit{a domain of a periodic type} of the period $q \geq 1$, if $f^{q} (S) = S$, $f^{k} (S) \cap S = \varnothing, k = 1... Q-1$.  If $q=1$ we have an invariant domain.
\end{intr}

\begin{intr}\label{intr18}
According to \cite{AM}, for a two-manifold the pathcomponents of a set $M \setminus (\textbf{cl} W') $ can be only of following types: \par
   1) simply-connected components of  wandering type \par
   2) simply-connected components of  periodic type \par
   3) doubly-connected components of  periodic type. \par 
\noindent
In the last case $M$ is $S^2$, $S$ is a component of  invariant type, and the boundary $\partial S$ consist of precisely two fixed points : $\partial S=\{\alpha, \omega \}$, where $\alpha$ -- source, $\omega$ -- sink, and $M=S \cup \{\alpha, \omega \} = S^{2}$. \par
\end{intr}

\begin{intr}
Let $\Omega_{k}$ be a periodic trajectory from $\Omega (f)$.
 According to \cite{ps}, the set of periodic trajectories  is partially ordered : $\Omega_{k} \leq \Omega_{l}$ if $W^{u} (\Omega_{l}) \cap W^{s} (\Omega_{k}) \not= \varnothing $. This partial order is also defined for periodic points. \par

\end{intr}

\begin{intr}
\textit{Definition} (1.7, \cite{Palis}). let $x, y \in \Omega ' (f) $ and $W^{u}(y) \cap W^{s}(x) \not= \varnothing $. Then, according to \cite{S1}, there exists the sequence of points $y_{0} = x, y_{1},.., y_{n} = y $ from $\Omega '(f)$ where $ y_{i + 1}$ does not belong to the trajectory of $ y_{i} $ such that $W^{s}(y_{i}) \cap W^{u} (y_{i + 1}) \not= \varnothing $. 

     Let us denote by $beh (y | x)$ the maximum length of such sequences. If $W^{u}(x) \cap W^{s}(y) = \varnothing $ then $ beh (y |x) = 0$. If we display the partial order on the set of periodic points with an oriented graph then $ beh (y | x)$ be the maximum length of oriented path from $y$ to $x$.
\end{intr}

\begin{intr}\label{intr17}
Let $p \in \Omega (f)$ be a fixed periodic point of $f$. According to \cite{ps}, a tubular neighborhood $T^{s}(p)$ of  the manifold $W^{s}(p)$ is a collection  $\{ t^{s}(\xi) \}$ of disjoint $C^{r}$-submanifolds $t^{s}(\xi)$ of $M$ indexed by $\xi$ from an open neighborhood $N$ of $p$ in $W^{u}(p)$ with the following properties : \par
    1) $V=V(W^{s}(p)) = \cup_{\xi \in N}t^{s}(\xi)$ is an open subset of $M$, containing $W^{s}(p)$; \par
    2) $ t^{s}(p) = W^{s}(p) $; \par
    3) $ t^{s}(\xi) $ intersects $N$ transversally  in the single point $\xi$; \par
    4) The map $\pi: V \rightarrow N$ which sends $t^{s}(\xi)$ in $\xi$ is continuous; the section $s$ which sends $\eta  \in  t^{s}(\xi)$ into the tangent space of $t^{s}(\xi)$ at $\eta$ is a continuous map  from $V$ into the Grassman bundle over $V$. 
\end{intr}

\begin{intr}
A tubular neighborhood $T^{s, k}$ of $W^{s}(\Omega_{k})$, where $ \Omega_{k}$ is a periodic trajectory of $f$, is the union of tubular neighborhoods $T(W^{s}(p)) = \{t^{s}(\xi)\}$ of a submanifold $W^{s}(p)$ for some $p \in \Omega_{k}$ and its images under the action of $f$. 
We denote the tubular neighborhood of $W^{s}(\Omega_{k})$ by $T^{s, k}$.
$$
T^{s, k}= \{f^{i} (t^{s}(\xi))\: |\: t^{s}(\xi)\in W^{s}(p),\: i = 0,\ldots, n-1\}\: \hbox{where $n$ is a  period of $\Omega_{k}$.}
$$
 The tubular set $T(W^{s}(p))=\{t^{s}(\xi) \}$ of a submanifold $W^{s}(p)$ is called invariant if $f^{- n}(t^{s}(\xi)) = t^{s}(f^{- n} (\xi)) $, where $n$ is a period of $p$. The tubular neighborhood $T^{s, k}$ of $W^{s}(\Omega_{k})$ is called invariant if it is obtained from an invariant tubular neighborhood of its point. A system of (stable) tubular neighborhoods of $f$ is the set of tubular neighborhoods constructed for  stable manifolds of each periodic trajectory. It is denoted by $\{T^{s, k}\}$. The same way defined is a system of unstable tubular neighborhoods of $f$ $\{T^{u, k}\}$. 
The system is called {\it invariant} if it consists of only invariant tubular neighborhoods.
The system is called {\it admissible} if it satisfies the following condition: from $t^{s}(\xi) \cap t^{s}(\eta) \not= \varnothing $ it follows that one submanifold contains another. \par

Palis and Smale have proved in \cite{Palis}, \cite{ps} that for $f$ there exist invariant admissible systems $\{T^{u, k}\}$ and $\{T^{s, k}\}$.
\end{intr}

\section{The basic principles of topological classification.}

Our first rather trivial observation is the fact that up to the topological conjugacy the information about Morse-Smale diffeomorphisms is concentrated on the set of non-regular and periodic points.
This set is the union of sinks, sources, and stable and unstable manifolds of saddle periodic points and is a one-dimensional complex. 

The complement to the set is the union of wandering and periodic components, where each of them is homeomorphic to a disk\footnote{Except the special case where $M=S^2$, see introduction, \ref{intr18}}. Every such component is determined by its edge, therefore if we consider a diffeomorphism up to the homeomorphism then, roughly speaking, we may consider this complex instead of the manifold and a permutation on the set of vertices and edges instead of the diffeomorphism, provided the subcycles corresponding to the wandering and periodic components are marked out. Note that a graph with marked out subcycles corresponding to the wandering and periodic components can be replaced with  the graph with spin (graph whose vertices are oriented circles), as it is done in this paper.

More strictly, the following theorem take place.

Consider the sets of periodic points, heteroclinic points, segments of stable and unstable manifolds of saddle periodic points between heteroclinic points, wandering and periodic components of a diffeomorphism $f$. Restrictions of $f$ on the corresponding sets are permutations on this sets.

\begin{theorem}\label{theoremf1} 
Let we have a mapping which maps periodic points, heteroclinic points, stable and unstable manifolds of saddle periodic points, wandering and periodic components of the diffeomorphism $f$ in corresponding objects of $g$, such that it establish a one-to-one correspondence between those objects of $f$ and $g$ and their permutations describing the action of the diffeomorphisms are conjugate with the help of this mapping. Then these diffeomorphisms are topologically conjugate.
\end{theorem}

The theorem is proved in the section \ref{Sec2Finit}.

We have that this one-dimensional complex indeed determine the diffeomorphism up to the topological conjugacy. 

This approach in the pure form was applied by A.~N.~Bezdenezhnykh to the gradient-like Morse-Smale diffeomorphisms in the paper \cite{Bgri1}. It become possible because there appeared a finite graph, but in the general case, which we consider, this one-dimensional complex has infinitely many vertices so another methods must be used.

\medskip
Our second observation is the fact that this one-dimensional complex possesses the local structure of direct product. As a matter of fact, the local structure of direct product appears on the set of heteroclinic and periodic points due to hyperbolicity, as one of the properties of locally maximum hyperbolic sets. 

In this paper we do not prove the fact that the set of heteroclinic and periodic points is a locally maximum hyperbolic one because the local structure of direct product follows directly from $\lambda$-lemma (see lemma \ref{lemmaf1}).

The local structure of direct product generates regular patterns in disposition of heteroclinic points. In this paper we use this patterns to specify position of heteroclinic points in the one-dimensional complex.

\begin{figure}[htb]
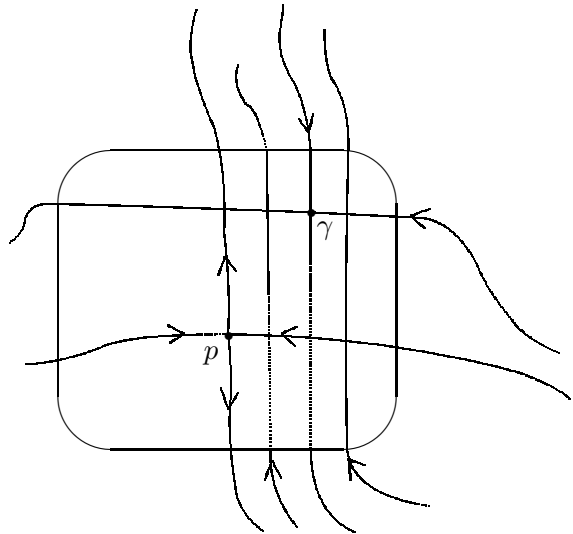

\begin{center}
\input FIG5.PIC
\end{center}
\caption{A neighborhood of $p$ with the structure of direct product.}\label{fig5}
\end{figure}

Our third observation gives us the way to use the structure of direct product in our construction. Consider a neighborhood of a periodic point $p$ with the structure of direct product (see fig. \ref{fig5}).  We see that in this neighborhood a heteroclinic point $\gamma$ and its neighbors on its  stable segment are located in the same order as their projections on stable manifold of $p$. Knowing the mutual location of projections on the stable and unstable manifolds of $p$, we can identically restore the rest of heteroclinic points in the neighborhood with the structure of direct product of $p$.

We see that internal points of the neighborhood are subordinate to their projections in the sense of the aforesaid statement. At that this projections also can have projections and so on. Our third observation is the fact that heteroclinic points which have no projection mold the rest of heteroclinic points and there is a finite number of such trajectories (see lemma \ref{lemmaLT1}).

We only need to remember periodic trajectories and heteroclinic trajectories which have no projections to be able to reconstruct the mutual location of the rest of heteroclinic trajectories using inductive process. At that we do not need to specify every vertice of our one-dimensional complex as a vertice of a graph, because the position of a heteroclinic point in the complex is uniquely determined by its position on its stable and unstable manifolds.
To specify the complex we separately write down the graph which specifies the incidence in periodic points and the set of heteroclinic points which are determined by their position on  stable and unstable manifolds of periodic points.

This is the general idea of the topological classification.

\subsection{Notes on numeric characteristics of heteroclinic points.}
To use reasoning on induction we must introduce a numeric characteristic of the subordination on the set of heteroclinic points. Since the numeric characteristic we use, the so called local type of heteroclinic points, is non-trivial, it would be better to say first a few words on the general ideas of such characteristics.
For this purpose we introduce the concept of lattice structure which is a good illustration of  numeric characteristics connected with the structure of direct product. On the other hand, this definitions will be used in section \ref{section3} where we rigorously define the local type of heteroclinic points.

First we define the set of subordinate heteroclinic points.

Let $p$ be a saddle periodic point of $f$ and $n$ be a period of $f$. Let $\gamma \in W^{u}(x)\cap W^{s}(y)$ be a heteroclinic point such that $W^{u}(x)$ intersects  $W^{s}(p)$ in a point $\gamma_1$ and $W^{s}(y)$ intersects $W^{u}(p)$ in a point $\gamma_2$.

\begin{figure}[htb]
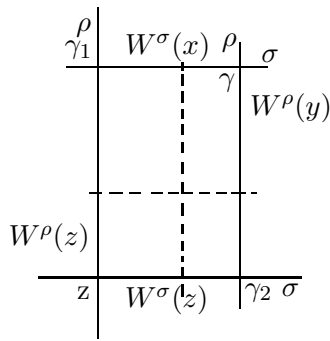

\begin{center}
\input FIG1S.PIC
\end{center}
\caption{A tetragon formed by heteroclinic point $\gamma$.}\label{fig1}
\end{figure}

\begin{definition}\label{definition200}
A heteroclinic point $\gamma \in W^{u}(x)\cap W^{s}(y)$, $p \not= x$, $p \not= y$ {\sl belongs to a lattice structure of $p$} if within the bounds of a tetragon formed by points $p$, $\gamma_{1}$, $\gamma$ and  $\gamma_{2}$ stable and unstable manifolds from $W'$ intersect one another only in one point.  
(see fig. \ref{fig1})
\end{definition}

In other words, the point $\gamma$ belongs to the neighborhood of $p$ with the local structure of direct product.

The point $\gamma_{2}$ is called a projection of $\gamma$ on $W^{\sigma}(y)$ and denoted by $\pi (\gamma, W^{\sigma}(y))$ or $\pi (\gamma)$ if the manifold is implied in context.

The tetragon formed by points $p$, $\gamma_{1}$, $\gamma$ and   $\gamma_{2}$ is called {\em the lattice tetragon} or {\em the tetragon with lattice structure} between $\gamma$ and $p$.

\begin{remark}
A heteroclinic point $\gamma$ can belong to lattice structures of different periodic points and can have more than one tetragon with the same periodic point.
\end{remark}

Such points always exist. 
For example,  heteroclinic points from the tubular neighborhood of $p$ with lattice structure belong to the lattice structure of $p$ according  to the definition.

\begin{lemma}\label{lemmaLT1}
There exists a finite number of heteroclinic trajectories that do not belong to the lattice structure of any periodic point.
\end{lemma}

\begin{proof}
Otherwise there exists infinitely many heteroclinic points with this property on a fundamental neighborhood and they must have a boundary point. At that all point in a small neighborhood of the boundary point must have projections.
\end{proof}

Using this definition, we can measure the degree of subordination of a periodic point $\gamma$ with the number of those periodic points who form with $\gamma$ tetragons with the lattice structure. 

But this numeric characteristic of subordination on the set of heteroclinic points require difficult and complicated inductive process. In this paper we propose a roundabout way.
In section \ref{section3} we introduce a more complicated numeric characteristic of subordination, a so called local type of heteroclinic points. First we define a so called lattice neighborhood of a periodic point $p$. It is the set of heteroclinic points which belong to the lattice structure of $p$ and satisfy some additional properties. We can imagine this set as one contained in an open neighborhood of $p$ containing the stable and unstable manifolds of $p$. Then we define a so called lattice neighborhood of a heteroclinic point $\gamma$. Unlike lattice neighborhoods of periodic points, a lattice neighborhood of a heteroclinic point is imagined as a set of heteroclinic points contained in a segment of stable or unstable manifold (see fig. \ref{fig6} on page \pageref{fig6}).

According to the sense of definition, all heteroclinic points which belong to the lattice neighborhood of a heteroclinic point $\gamma$ are subordinate to $\gamma$. At that for the lattice neighborhoods of heteroclinic points the following property take place (see remark \ref{remark28}): if two lattice neighborhoods of heteroclinic points on the same manifold have non-empty intersection then one of them contains the other. It considerably simplifies the mutual subordination on the set of heteroclinic points: if a point $\gamma_1$ is subordinate to a point $\gamma_2$ and the point $\gamma_2$ is subordinate to a point $\gamma_3$ then the point $\gamma_1$ is also subordinate to the point $\gamma_3$.

We define local type of a point to measure a number of lattice neighborhoods which contain the point. At that so defined local type has a very useful property (lemma \ref{lemma36plus}): the local type of a heteroclinic point is equal to sum of local types of its projections in any lattice structure to which the point belongs.

Denote the local type of a heteroclinic point $\gamma$ by $loctype(\gamma)$. Let $\gamma$ belongs to lattice neighborhood of a periodic point $p$ and $\gamma_1$ and $\gamma_2$ be projections of $\gamma$ in lattice structure of $p$. We have
\begin{equation}
loctype(\gamma)=loctype(\gamma_1)+loctype(\gamma_2)
\end{equation}
This equation considerably simplifies the further construction.

\section{Numbering of heteroclinic points and coding sets for heteroclinic trajectories.}
\label{Sec4defCS}

In this section we construct identifiers for heteroclinic points. Every identifier (so called simple formula) uniquely determine the relative position of the heteroclinic point on stable and unstable manifolds and a direction of the intersection of stable and unstable manifolds in this heteroclinic point. 

\subsection {Numbering of heteroclinic points.} 
To describe a heteroclinic point we first assign to it an index number.
We do it in the following way: 
 On every $W_{i}^{\sigma}(x_{j})$ we first arbitrarily choose a start numeration point of local type 1 and denote it by $P_{ij}^{\sigma}(0)$.  (Here $P_{ij}^{\sigma}$ is called a prefix, and $(0)$ is called a vector of the point number). Let  $m$ be a minimum number such that a manifold $W_{i}^{\sigma}(x)$ is invariant under the action of $f$. A fundamental neighborhood [$P_{ij}^{\sigma}(0), f^{m}(P_{ij}^{\sigma}(0))] \subset W_{i}^{\sigma}(x_{j})$
contains $N_{ij}^{\sigma} + 1$ points of local type 1 (see \ref{lemma34}). Let us denote them by $P_{ij}^{\sigma}(0)$, $P_{ij}^{\sigma}(1)$,\dots, $P_{ij}^{\sigma}(N_{ij}^{\sigma})$, numbering  ``on current ''. Let points of the local type 2 are contained in the lattice neighborhood of a point $P_{ij}^{\sigma}(l_{0})$. According to (\ref{lemma37}), their set is infinite and $P_{ij}^{\sigma}(l_{0})$ is a unique boundary point of this set. Then the point $P_{ij}^{\sigma}(l_{0})$ divides their set into two subsets (one of them can be empty): a front (''on current``) set and a back set. Let us choose in the front subset the extreme point and denote it by $P_{ij}^{\sigma}(l_{0}, -1)$. Other points of the front set we number ''against current'' from a point $P_{ij}^{\sigma}(l_{0}, -1)$ in the direction of $ P_{ij}^{\sigma} (l_{0}):$
$$
 P_{ij}^{\sigma} (l_{0}, -1), P_{ij}^{\sigma} (l_{0}, -2), P_{ij}^{\sigma}(l_{0}, -3), \dots
$$
Let us remark that in this case we already have a two-component vector. 

  The same we do with the back set. We denote extreme point by $P_{ij}^{\sigma}(l_{0},1)$ and number other points of the back set ''on current'' from the point $P_{ij}^{\sigma} (l_{0},1)$ in the direction of $P_{ij}^{\sigma}(l_{0}):$
$$
P_{ij}^{\sigma}(l_{0},1),P_{ij}^{\sigma}(l_{0},2), P_{ij}^{\sigma}(l_{0},3), \dots
$$

The further construction is inductive. If  a lattice neighborhood of already numbered point $P_{ij}^{\sigma}(l_{0}, l_{1}.., l_{n})$ of the local type $n$ contains points of the local type $n+1$,  we number points of the front set exactly in the same way ''against current'': 
$$
P_{ij}^{\sigma}(l_{0}, l_{1}.., l_{n}, -1), P_{ij}^{\sigma}(_ {} l_{0}, l_{1}.., l_{n}, -2),  P_{ij}^{\sigma}(_ {} l_{0}, l_{1}.., l_{n}, -2), \dots
$$
And we number points of the back set ''on current'':
$$
P_{ij}^{\sigma}(l_{0}, l_{1}.., l_{n}, 1), P_{ij}^{\sigma} (l_{0}, l_{1}.., l_{n}, 2),  P_{ij}^{\sigma} (l_{0}, l_{1}.., l_{n}, 3), \dots
$$

  Then the numbers constructed in this way  we spread outside the fundamental neighborhood. We assign numbers to points, using already constructed numbers of their images under the action of  $f^{m}$ in the following way:  if ${(f^{-m})}^{k}(\gamma) = P_{i, j}^{\sigma}(l_{1}, l_{2}.., l_{n})$ then we denote $\gamma$ by $P_{i, j}^{\sigma}(l_{1} + k*N_{ij}^{\sigma}, l_{2}.., l_{n})$. Thus, every heteroclinic point obtains a number in a form of a sequence of integers, called a vector of point numbers or simply vector. 
The lexicographic order on the set of vectors of points with an identical prefix coincides with an order of points on a submanifold corresponding to the prefix. 

 Since the amount of heteroclinic points in lattice neighborhoods is a topological invariant, all components of the point's vector except the first one are its topological invariants.  The begin numeration point  $P_{ij}^{\sigma}(l_{0})$ is chosen arbitrarily, therefore on each $W_{i}^{\sigma}(x_{j})$ the first components of point number vector are defined up to a shift on an arbitrary integer $m^{\sigma}_{ij}$, identical to all heteroclinic points on $W_{i}^{\sigma}(x_{j})$. \par
\subsection {Formulas and coding sets.} 
Now we are ready to describe a heteroclinic point.

  Let us associate with a heteroclinic point $\gamma$ from $W_{i}^{u}(x_{k}) \cap W_{j}^{s}(x_{l})$ a three-component set
$(P^{u}_{ik}(t_{0}, t_{1}..., t_{n}), P^{s}_{jl}(g_{0}, g_{1}..., g_{n}), D)$, where $n$ is local type of the point, $D$ accepts the value $+1$ if the orientation of a frame formed  by manifolds $W_{i}^{u}(x_{k})$ and $W_{j}^{s}(x_{l})$ in a point $\gamma$ coincides with the orientation of the manifold, and $D$ accepts the value $-1$ if the orientations do not coincide. This set is called \textit{a simple formula of the heteroclinic point}. \par

\begin{definition}\label{definition52}
A formula of the form 
$$
(P^{u}_{ik}(t_{0}+ k* N_{ik}^{u}, t_{1},\ldots, t_{n}), P^{s}_{jl} (g_{0} + k* N_{jl}^{s}, g_{1},\ldots, g_{n}), D)_{k \in \zz}, 
$$
 obtained from a  point $(P^{u}_{ik} (t_{0}, t_{1},\ldots, t_{n}), P^{s}_{jl} (g_{0}, g_{1},\ldots, g_{n}), D)$ with adding a displacement under the action of the diffeomorphism,  is called \textit {a simple formula of a heteroclinic trajectory}. 
\end{definition}

\begin{definition}\label{definition51}
{\em
The next expression is called  {\em a formula of a heteroclinic point set} or, simply, formula. 
\begin{displaymath}
\left(
\begin{array}{ccc}
P^{u}_{ik}\bigl(t_{0} (p_{0},\dots, p_{r}),\dots, t_{n} (p_{0},\dots, p_{r})\bigr), \\ 
P^{s}_{j1} \bigl(g_{0} (p_{0},\dots, p_{r}),\dots, g_{n} (p_{0},\dots, p_{r})\bigr), \\
 D (p_{0},\dots, p_{r})
\end{array}
\right)
\end{displaymath}
 Here $t_{i}$, $g_{i}$, $D$ are linear functions  of  parameters $p_{0},\dots, p_{r}$, where $p_{0} \in \zz$ is called a trajectory parameter, 
$p_{1},\dots, p_{r}$ is called local parameters and the domain of definition of $p_{1},\dots, p_{r}$ is a subset of $\zz^r$ determined by a finite set of linear inequalities.
}\end{definition}

   The formula associates with each set of numbers $(p_{ },\dots, p_{r})$ some heteroclinic point from $W_{i}^{u}(x_{k}) \cap W_{j}^{s}(x_{l}) $. \par

\begin{definition}\label{definition55}
a set of simple formulas which contains simple formulas of every trajectory of local type 1 and simple formulas of extreme trajectories of its lattice neighborhoods  for this trajectory is called \textit {a basic coding set}.
\end{definition}
In other words, it contain four formulas of extreme trajectories  of local type 2 in back and front sets in stable and unstable manifolds of this trajectory if they are not empty, four formulas of local type 3 if they are not empty and so on up to the maximal  local type.

According to the definition,  the basic coding set is finite.

\begin{definition}\label{definition53}
 A set of formulas describing all trajectories of the diffeomorphism, is called \textit {the extended coding set}. \par
\end{definition}

\begin{definition}\label{definition54}
{\em
 Two coding sets is called isomorphic if there exist:

i) a modification of the manifold orientation (we denote it $r \mycirc D$ assuming that $r \mycirc D=+ D$ if orientation does not change and  $ r \mycirc D=-D$ otherwise)

ii) the set of automorphisms of numeration for every $W_{i}^{\sigma}(x_{j})$, which are determined by integers $n^{\sigma}_{ij}$

iii) a renumbering  $\varphi$ of the set of periodic points of $f$ and a renumbering  $\psi$ of the set of manifolds from $W ' \setminus \Omega (f)$

such that if the first coding set determines a point 
$$
(P^{u}_{ik} (t_{0}, t_{1},..., t_{n}), P^{s}_{jl} (g_{0}, g_{1}..., g_{n}), D)
.$$
 then the second extended coding set determines a point 
$$
(P^{u}_{\varphi (i), \psi (k)} (t_{0} + n^{u}_{ik}, t_{1}. ., t_{n}), P^{s}_{\varphi (i), \psi (l)} (g_{0} + n^{s}_{jl}, g_{1}..., g_{n}), r (D)) 
$$
 and vice versa. \par
}\end{definition}

The basic and  extended coding sets are the topological invariants of the diffeomorphism up to their isomorphisms. 

\subsection{Finiteness of the extended coding set.}

The main part of the information we store in the topological invariant describes a mutual location of periodic and heteroclinic points of the diffeomorphism. We describe it using formulas defined above.

There are infinitely many heteroclinic points on the manifold but due to hyperbolic structure of the diffeomorphism there appears an order in their disposal which we call a lattice structure. This order is essential only in small neighborhoods of stable and unstable manifolds of  every saddle periodic point, far off them it gets broken. But system of those neighborhoods contains all heteroclinic points. 

As proved in lemma \ref{lemma34} there is a finite number of trajectories of local type 1. They are ``nodes of fracture'' of lattice structure of diffeomorphism. Positions of the rest of points entirely repeat the order of position of their projections in lattice neighborhoods owing to lattice structure. Since projections have a smaller local type then origin (lemmas \ref{lemma35} and \ref{lemma36}), we are able to describe all heteroclinic trajectories on induction, starting from trajectories of local type 1 and trajectories that determine the shape of lattice neighborhood,  stored in the basic coding set.   
    This reasons suggest that the extended coding set can be written down as a finite number of the formulas, but does not give an algorithm in an explicit aspect. The next theorem gives us a constructive proof of this fact.

\begin{theorem}\label{TheoremCS2}
The basic coding set and the graph $G(f)$\footnote{see next section} (both finite) are uniquely determine the extended coding set. The obtained extended coding set is also finite.
\end{theorem}

The proof of the theorem is given in section \ref{Sec5AlgCS}. The proof gives us the algorithm of construction of the extended coding set on the basis of the basic coding set and the graph defined below. Obtained extended coding set of formulas is finite.
\section{Distinguishing graph of the diffeomorphism.} 
\label{Sec7graph}

Let us associate with the diffeomorphism $f$ an oriented graph with spin $G(f)$ (a graph where every vertice is an oriented circle), whose vertices correspond to points from $\Omega(f)$ and oriented edges correspond to oriented ``on current'' components of connectivity of the set $W' \setminus \Omega'(f)$.  All vertices are oriented in correspondence with the orientation of the manifold.  Then the graph has the following properties:

\begin{enumerate}
\item the vertices of the graph for which all  incident edges  go in correspond to sinks;
\item the vertices for which all incident edges go out correspond to sources;
\item the remaining vertices correspond to saddles.
\item every saddle vertice is exactly incident to a pair of going in (stable) edges and to a pair of going out (unstable) edges  located crosswise on  an oriented circle which is the vertice.
\item Each edge of the graph is always incident by one of its ends to some saddle vertice and  is incident by another end either to a sink, or to a source, or is free with this end (this case corresponds to a case when the corresponding manifold contains heteroclinic points. The corresponding edge is called {\em a free edge}). 
\end{enumerate}

We assign to a vertice the weight $``+''$, if its orientation type (see \ref{intr12}) is equal to $1$, and  the weight $``-''$ otherwise.
The diffeomorphism $f$ induces a automorphism $ S_{f}:G (f) \rightarrow G(f)$ on the set of vertices of the graph. It is described by a permutation on the set of vertices. 

\begin{definition}\label{definition81}
The graph $G(f)$ with the corresponding basic coding set and permutation $S_f$  where edges corresponding to manifolds containing heteroclinic points and incident to them saddle vertices are marked by their numbers in the basic coding set is called \textit{a distinguishing graph  of $f$} and is denoted by  $G^*(f)$.
\end{definition}

\begin{definition}\label{definition}
Distinguishing graphs $G^*(f)$ and $G^*(g)$ of diffeomorphisms $f$ and $g$ respectively are isomorphic, if there exist
\begin{enumerate}
\item isomorphism $\varphi$ of graphs $G(f)$ and $G(g)$, preserving weights and orientation of edges and vertices, such that $S{f} = \varphi  S_{g} \varphi^{-1}$;
\item a change of orientation $r$, which either preserves orientations of vertices of the graph and orientation components of formulas in the basic coding set or simultaneously reverse the orientation of vertices of the graph and orientation components of formulas in the basic coding set; 
\item a set of automorphisms of numeration for every free edge of $G(f)$
\end{enumerate}
such that basic coding sets for $f$ and $g$ are isomorphic. 
\end{definition}

    Thus, we have associated with a diffeomorphism the distinguishing graph which is a finite object. Isomorphism of two distinguishing graphs specify the one-to-one correspondence between periodic points, heteroclinic points, stable and unstable to manifolds of saddle periodic points, wandering and periodic components of the diffeomorphisms. According to the theorem \ref{theoremf1}, it leads to the topological conjugacy of $f$ and $g$ and the following theorem takes place.

\begin{theorem}\label{theorem3Main}
Diffeomorphisms $f$ and $g$ are topologically conjugate if and only if their distinguishing graphs are isomorphic.
\end{theorem}

\section {Example of the invariant.}

\begin{figure}[htbp]
\unitlength=1mm
\includegraphics [width=11cm, keepaspectratio]{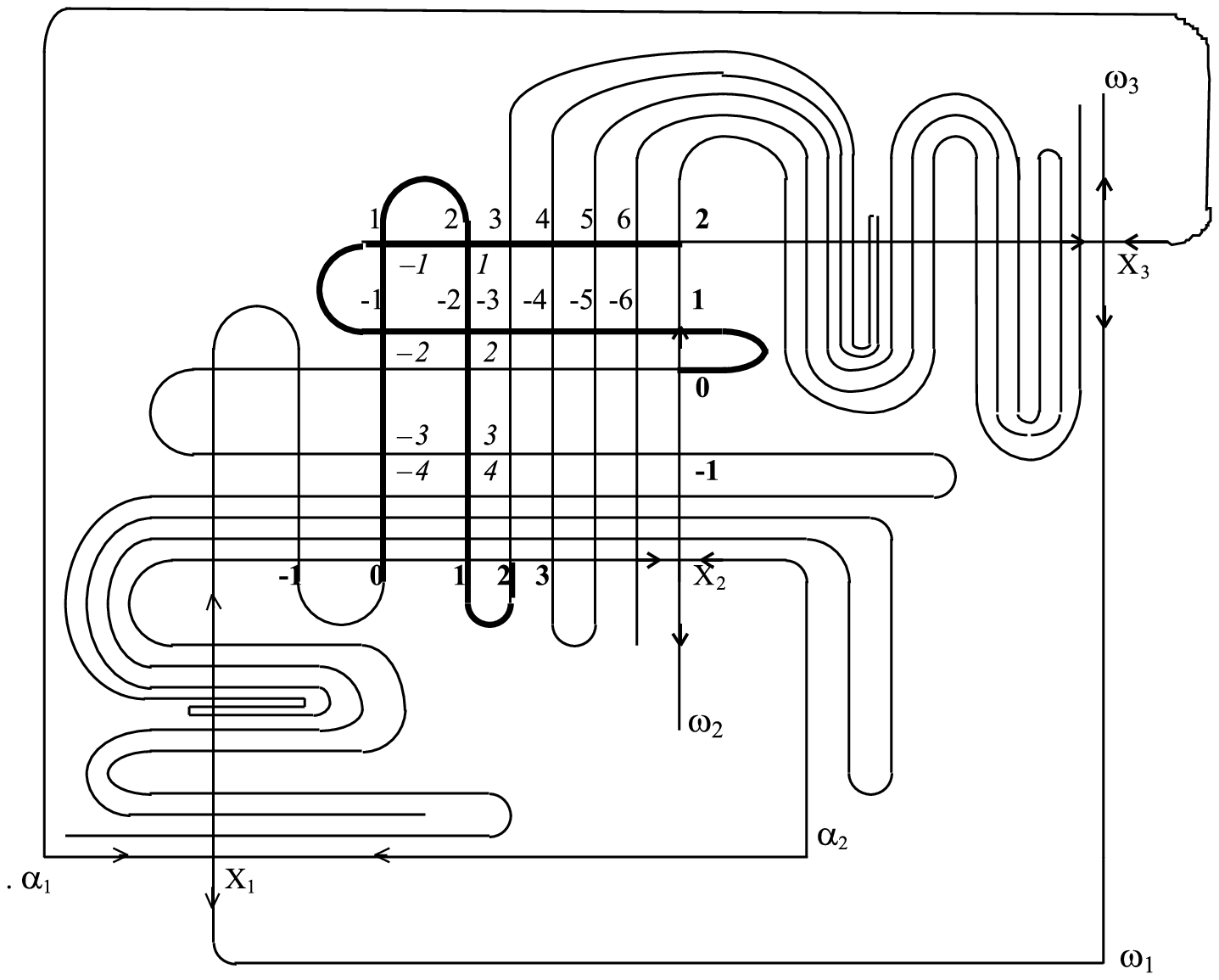}
%
\hbox{
\vbox{\noindent
\divide \hsize by 2                                                                             

The basic coding set is: \\
$(P^{u}_{11}(2n),P^{s}_{12}(2n),+1)_{n\in \zz};$\\
$(P^{u}_{11}(2n+1),P^{s}_{12}(2n+1),-1)_{n\in \zz};$\\
$(P^{u}_{12}(2n),P^{s}_{13}(2n),+1)_{n\in \zz};$\\
$(P^{u}_{12}(2n+1),P^{s}_{13}(2n+1),-1)_{n\in \zz};$\\
$(P^{u}_{11}(2n,-1),P^{s}_{13}(2n+2,1),+1)_{n\in \zz};$\\
$(P^{u}_{11}(2n+1,1),P^{s}_{13}(2n+2,2),-1)_{n\in \zz};$\\
$(P^{u}_{11}(2n,-2),P^{s}_{13}(2n+1,-1),+1)_{n\in \zz};$\\
\\The permutation $S_f$ is trivial. \\
}
\input FIG4.PIC
} 
\caption{The Morse-Smale diffeomorphism of a sphere.
 Its distinguishing graph and basic coding set}
\label{fig3}
\end{figure}

\begin{figure}[htb]
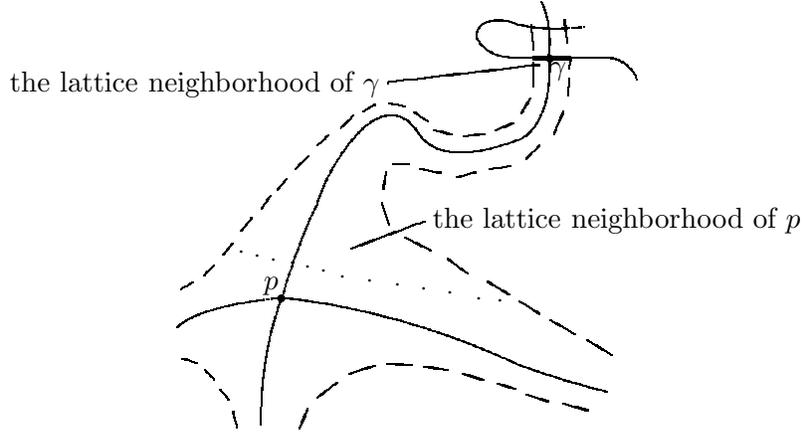

\begin{center}
\input FIG6.PIC
\end{center}
\caption{The lattice  neighborhoods of $p$ and $\gamma$.}\label{fig6}
\end{figure}

As an example we consider a Morse-Smale diffeomorphism of a sphere which have with two fixed sources, three fixed sinks, three fixed saddle points and heteroclinic trajectories of the local types 1 and 2. All periodic points have the orientation type 1. This example is illustrated by fig. \ref{fig3}. Fundamental neighborhoods are emphasized by bold lines. The numeration of heteroclinic points of local type  1 is in bold style, the  numeration of heteroclinic points of local type 2 on $W^u_1(x_1)$ is in italic style, the  numeration of heteroclinic points of local type 2 on $W^u_1(x_3)$ is in plain style.

\medskip \noindent
The extended coding set of the diffeomorphism on fig. \ref{fig3} is: \\
$(P^{u}_{11}(2n),P^{s}_{12}(2n),+1)_{n\in \zz}; (P^{u}_{11}(2n+1),P^{s}_{12}(2n+1),-1)_{n\in \zz};$\\
$ (P^{u}_{12}(2n),P^{s}_{13}(2n),+1)_{n\in \zz}; (P^{u}_{12}(2n+1),P^{s}_{13}(2n+1),-1)_{n\in \zz};$\\
$(P^{u}_{11}(2n+2k,-1-2k),P^{s}_{13}(2n+2,2k+1),+1)_{n\in \zz, k\geq 0};$\\
$(P^{u}_{11}(2n+2k+1,1+2k),P^{s}_{13}(2n+2,2k+2),-1)_{n\in \zz, k\geq 0};$\\
$(P^{u}_{11}(2n+2k,-2-2k),P^{s}_{13}(2n+1,-2k-1),+1)_{n\in \zz, k\geq 0};$\\
$(P^{u}_{11}(2n+2k+1,1+2k),P^{s}_{13}(2n+1,-2k-2),-1)_{n\in \zz, k\geq 0};$\\

\section{Topological conjugacy of Morse-Smale diffeomorphisms.}
\label{Sec2Finit}

This section is devoted to the proof of the theorem \ref{theoremf1}.

The proof uses methods similar to methods of J. Palis and S. Smale's paper \cite{ps} where they used tubular neighborhoods (see \ref{intr17}) to construct a homeomorphism which conjugate Morse-Smale diffeomorphisms.

To simplify the proof, we use tubular neighborhoods with imposed additional conditions.

\begin{definition}
A systems $\{T^{u, k}\}$ and $\{T^{s, k}\}$ of tubular neighborhoods is called systems of tubular neighborhoods with lattice structure if they are
1. {invariant;} 

2. {admissible;} 

3. {all fibers of the systems are connected;} 

4. {every tubular neighborhood $T^{u, k}$ coincides with the corresponding $T^{s, k}$ as the topological neighborhood;} 

5. {if $t^{u, k}(\xi)\cap t^{s, l}(\eta) \not= \varnothing $ then they have a unique point of intersection.} 

\end{definition}

One readily sees that the union and the intersection of the systems of tubular neighborhoods with lattice structure again is a system with lattice structure.

The existence of systems of tubular neighborhoods with this properties was implicitly mentioned in J. Palis and S. Smale's paper \cite{ps}. A simple proof of the fact is given here also.

\begin{lemma}\label{lemmaf1}
Every Morse-Smale diffeomorphism on two-dimensional closed manifold has a system of tubular neighborhoods with lattice structure.
\end{lemma}

\begin{proof}
Let $ \Omega_{k} \in \Omega'(f)$ be a periodic saddle trajectory and $ T^{u, k}, T^{s, k}$ be the tubular neighborhoods of its stable and unstable manifolds and $x$ be a periodic point of $\Omega_{k}$. 
According to the $\lambda $-lemma for families of transversal disks (see \cite{PalisMelo}), $\forall \varepsilon > 0 $ there exists a small neighborhood of the point $x$ where the fibers from $T^{u, k}$  are $\varepsilon$-$C^{1}$ -- close to $W^{u}(x_{k}) $ and the fibers from $T^{s, k}$ are $\varepsilon$-$C^{1}$ --  close to $W^{s}(x_{k})$. 
Choosing a small $ \varepsilon $, we obtain a small neighborhood $U(x)$ of the point $x$, $U(x)\subset T^{u, k} \cap T^{s, k}$, such that the fibers from $ T^{u, k} $ and $ T^{s, k} $ intersect each other transversally in the unique point. 

Choose on $ W^{u}(\Omega_{k})$ and $W^{s}(\Omega_{k})$ fundamental neighborhoods $G^{u}(\Omega_{k})$ and $ G^{s}(\Omega_{k})$ such that the set $Q = \{x\in t^{s}(\xi) \cap t^{u}(\eta)\:|\: \xi \in G^{u}(\Omega_{k}) \eta \in G^{s}(\Omega_{k}) \}$ is contained in $U(x)$. $G^{u}(\Omega_{k})$, $G^{s}(\Omega_{k})$ consist of two pathcomponents  $G_{1}^{\sigma}(x)$ and $G_{2}^{\sigma}(x)$. Consider arbitrarily chosen  homeomorphisms $d_{i,j}: G_{i}^{u}(x) \rightarrow G_{j}^{s}(x)$ reversing the natural orientation of $G_{i}^{u}(x)$ and  $G_{j}^{s}(x)$. Let us extend $d_{i,j}$ to a homeomorphism $d_{i,j}: W_{1}^{u}(\Omega_{k}) \rightarrow W_{1}^{s}(\Omega_{k})$, putting $d_{i,j}(f^{n}(x)) = f^{n}(d_{i,j}(x))$, where $x\in G_{i}^{\sigma}(x)$, $n\in \zz $. 
Let $D$ be set of the graphs of those diffeomorphisms in $Q$, 
$$
i.e. \quad D= \{x\in t^{s}(\xi) \cap t^{u}(\eta)\ |\ \xi \in W^{u}(\Omega_{k})\ \eta \in W^{s}(\Omega_{k})\ d_{i,j} (\xi) = \eta \}.
$$ 
Let the new tubular neighborhoods with lattice structure $T_{1}^{k, u}$ and $T_{1}^{k, s}$ be pathcomponents of the sets $T^{k, u} \setminus D$ and $T^{k, s} \setminus D$ containing points from $\Omega_{k}$. 
The definition is satisfied because the properties of tubular neighborhoods hold out due to inheritance.
Repeating this procedure for other manifolds we obtain a system of tubular neighborhoods with lattice structure. 
\end{proof}

Further we will use terms ``tubular neighborhood'' and ``system of tubular neighborhoods'' implying that they are systems of tubular neighborhoods with lattice structure and tubular neighborhoods with lattice structure.

\medskip
The definition of $behavior$ of periodic points given in \cite{Palis}, \cite{ps} \ being introduced to describe periodic points requires mentioning of manifolds explicitly. It is convenient to use the definition describing heteroclinic points in the following form:
.
\begin{definition}\label{definition3}
{\em 
A point of intersection $W^{u}(y)$ and $W^{s}(x)$, $x, y\in \Omega'(f)$ is called {\em a point  of $beh$-type $n$} if $beh(y | x) = n$. The trajectory of the point is called {\em the trajectory of $beh$-type n}
}\end{definition}

\begin{definition}\label{definition4}
{\em
Let a heteroclinic point $\gamma$ belongs to the intersection $W^{\sigma}(z)$ and $W^{\rho}(x)$, $x, z\in \Omega ' (f) $. {\em A tubular neighborhood $UT^{\rho}(\gamma)$ of the heteroclinic point $\gamma$ in $W^{\sigma}(z)$} is a neighborhood of a point $\gamma$ in $W^{\sigma}(z)$ which is a pathcomponent of the set $T^{\rho}(W^{\rho}(x)) \cap W^{\sigma}(z)$ containing a point $\gamma $.  {\em A tubular neighborhood $UT^{\rho}(\Gamma)$ of the trajectory $\Gamma(\gamma)$} is the union of tubular neighborhoods of all points that belongs to $\Gamma$.
}\end{definition}

\begin{lemma}\label{lemmaf5}
Let $W^{\rho}(x)$ intersects $W^{\sigma}(y)$ with the $beh$-type $n$ (i. e. either $beh(x|y)$ or 
$beh(y|x)$ is equal to $n$). Then $W^{\rho}(x)$ cannot have a boundary curve, intersecting $W^{\sigma}(y)$ with a $beh$-type, large or equal $n$.
\end{lemma}

\begin{proof}
Suppose that there exists a boundary for $W^{\rho}(x)$ curve $W^{\rho}(z)$ which intersects $W^{\sigma}(y)$ with a $beh$-type, large or equal $n$. Then we can connect $x$ and $y$ with a way of length at least n. Since $W^{\rho}(z)$ is a boundary curve for $W^{\rho}(x)$ we can connect $x$ and $z$ with a way of length at least 1. Hence, we can connect $x$ and $y$ with a way of length greater or equal $n+1$. The maximum length of the way according to definition is $n$, we obtain contradiction.
\end{proof}

\begin{lemma}\label{lemmaf6}
For every saddle periodic point $y$ we always can choose a finite number of trajectories on its manifold $W^{\sigma}(y)$ such that they contain in their tubular neighborhoods in $W^{\sigma}(y)$ the rest of heteroclinic trajectories on $W^{\sigma}(y)$.
\end{lemma}

\begin{proof}
Consider the heteroclinic points of $beh$-type 1. $W^{\sigma}(y)$ contains a finite number of trajectories of $beh$-type 1. In other words, we can choose a finite number of trajectories  on $W^{\sigma}(y)$ of $beh$-type, not exceeding $1$, which contain in their tubular neighborhoods in $W^{\sigma}(y)$ the rest of heteroclinic trajectories on $W^{\sigma}(y)$ of $beh$-type, not exceeding $1$. 

We proceed by induction on $beh$-type of points on $W^{\sigma}(y)$. 
Suppose we can choose a set $\{z_{i}\}_{i = 1.. p}$ which consists of a finite number of trajectories  on $W^{\sigma}(y)$ of $beh$-type, not exceeding $k-1$, which contain in their tubular neighborhoods in $W^{\sigma}(y)$ the rest of heteroclinic trajectories on $W^{\sigma}(y)$ of $beh$-type, not exceeding $k-1$.
According to lemma \ref{lemmaf5} the set of points of $beh$-type $k$ cannot have boundary points of larger $beh$-type, therefore the set of trajectories of $beh$-type $k$ outside of tubular neighborhoods of trajectories of the set $\{z_{i}\} _{i = 1..p}$ is finite and can be described as $\{z_{i}\}_{i = p+1.. q}$.
Thus, the set $\{z_{i}\} _{i = 1..q}$ has the property that the tubular neighborhoods of its trajectories contain the rest of heteroclinic trajectories on $W^{\sigma}(y)$ of $beh$-type, not exceeding $k$.
We obtain the statement of the lemma when $k$ is the maximal  $beh$-type of points on  $W^{\sigma}(y)$.
\end{proof}

\begin{vtheorem}{Theorem}{\ref{theoremf1}} 
Let we have a one-to-one correspondence between periodic points, heteroclinic points, stable and unstable manifolds of saddle periodic points, wandering and periodic components of two diffeomorphisms $f$ and $g$ and their permutations describing the action of the diffeomorphisms are conjugate. Then these diffeomorphisms are topologically conjugate.
\end{vtheorem}

\begin{proof}
the main goal of the proof is to construct a homeomorphism that conjugate $f$ and $g$. In a special case when this diffeomorphisms have a finite number of heteroclinic trajectories such a homeomorphism was constructed by A. Z. Grines in the paper \cite{Gri1}. We consider a general case of infinite number of heteroclinic trajectories. According to the condition of the theorem, there exist a mapping $h_{g}\colon (W'{}^{u}(f) \cap W '{}^{s}(f)) \cup \Omega (f) \rightarrow (W '{}^{u}(g) \cap W '{}^{s}(g)) \cup \Omega (g)$, which translates the sets of periodic and heteroclinic points of $f$ in the corresponding sets for $g$. During the proof we will extend domain of definition of this mapping on the whole manifold. 

\medskip
\noindent
{\bf Part 1. }
In this part we construct mapping $h_{w}$ that extends $h_{g}$ on the one-dimensional complex $W ' (f) \cup \Omega (f)$. 

To simplify designations,  we put $h_{g}(\Omega _{k})=\Omega _{k}$.
We mark with tilde objects corresponding to $g$. Objects corresponding to $f$ we leave unmarked.

We construct $h_{w}$ using induction on size of the maximal  $beh$-type of points on $W^{\sigma}(\Omega_{k})$. 
A step $k$ of induction consist of two stages. On the first stage we extend $h_{g}$ on all manifolds having the maximal  $beh$-type of points equal to $k$. Then, on the second stage, we narrow or ``trim'' the tubular neighborhoods of heteroclinic points, cutting off their unnecessary parts. The trimmed tubular neighborhood in $W'{}^{\sigma}$ of a heteroclinic point $\gamma \in W^{\rho}(\Omega_{k})\cap W'{}^{\sigma}$ we denote by $T_{1}U^{\rho, k}(\gamma)$. 

We need the last stage in the step of induction because there is an ambiguity in size of the tubular  neighborhoods, so chosen systems of tubular neighborhoods for $f$ and $g$ may be different, though in any case they possess common parts. To avoid this difficulty, we construct during the part 1 of the proof systems of new tubular neighborhoods of heteroclinic points for $f$ and $g$ on the basis of common parts of old ones, such that the new tubular neighborhoods are in one-to-one correspondence. Then in the part 2 of the proof using obtained new tubular neighborhoods of heteroclinic points (they are sets of segments) we construct systems of new tubular neighborhoods for periodic points of $f$ and $g$ which also are in one-to-one correspondence (they are sets of open neighborhoods) and  use them to expand the homeomorphism on the whole surface.

Denote by $\nu_{\gamma}^{-1}$ a continuous mapping $TU^{\rho}(\gamma) \rightarrow \pi (TU^{\rho}(\gamma))$, which is the restriction of  the projection $\pi \colon TU^{\rho}(W^{\rho} (\Omega_{k})) \rightarrow W^{\sigma}(\Omega_{k})$ on points of  the segment $TU^{\rho}(\gamma)$. Since $TU^{\rho}(\gamma)$ is obtained as the intersection of $TU^{\rho}(\Omega_{k})$ with $W^{\sigma}(x)$ then $W^{\sigma}(x)$ intersects each fiber from $TU^{\rho}(\Omega_{k})$ in a unique point, so $\nu_{\gamma}^{-1}$ is a bijection and it has an inverse mapping  $\nu$.

\medskip
\noindent
{\bf Step 0 (The base of induction). }
Consider $W^{\sigma}(\Omega_{k})$ which contain no heteroclinic points (so the maximal  $beh$-type is equal to 0). In this case we simply take any fundamental neighborhoods $G^{\sigma}(\Omega_{k}) $ and $\widetilde G^{\sigma}(h_{g}( \Omega_{k}))$ and arbitrary homeomorphism $h_{0}: G^{\sigma}(\Omega_{k}) \rightarrow \widetilde G^{\sigma}(h_{g}( \Omega_{k}))$ that translates pathcomponents of $G^{\sigma}(\Omega_{k})$ in corresponding ones of $\widetilde G^{\sigma}(h_{g}( \Omega_{k}))$, preserving their orientation with respect to $\Omega_{k}$. 
Define $h_{w}$ on $W^{\sigma}(\Omega_{k})$ with a relation $h_{w}(x) = g^{n}(h_{g}(f^{-n}(x)))$ where $n$ is such that $f^{n}(x)\in G^{\sigma}(\Omega_{k})$. 

On the second stage we trim tubular neighborhoods.
The manifold $W^{\rho}(\Omega_{k})$ dual to $W^{\sigma}(\Omega_{k})$ may have  heteroclinic trajectories. To trim, we simply put  
$$
T_{1}U^{\rho, k}(\gamma) = \nu ( \nu^{-1}(TU^{\rho, k}(\gamma)) \cap h_{w}^{-1} ( \widetilde \nu^{1}(\widetilde {TU}^{\rho, k}(h_{g}( \gamma))))) 
$$
for each heteroclinic point $ \gamma \in W^{\rho}(\Omega_{k})$ and, similarly, 
$$ 
\widetilde {T_{1}U}^{\rho, k}(\gamma) =\widetilde \nu (h_{w}(\nu^{-1}(TU^{\rho, k}(\gamma))) \cap \widetilde \nu^{-1}(\widetilde {TU}^{\rho, k}(h_{g}( \gamma))) )
$$
for each heteroclinic point $ \gamma \in \widetilde W^{\rho}(\Omega_{k})$. $h_{w}$ is used correctly because it is already defined on $W^{\sigma}(\Omega_{k})$ in the previous stage.
This relations mean that we have cut off points which have either no image or no pre-image for the mapping $\widetilde \nu \circ h_{w}\circ \nu^{-1}$. After cutting off we obtain new tubular neighborhoods $T_{1}U^{\rho, k}(\gamma)$ and $\widetilde {T_{1}U}^{\rho, k}(h_{g}(\gamma))$ which are homeomorphic now with the homeomorphism $\widetilde \nu \circ h_{w}\circ \nu^{-1}$

Those operations of stages 1 and 2 we repeat for all manifolds which contain no heteroclinic points.

\medskip
\noindent
{\bf A step of induction.}
Suppose we have extended $h_{g}$ to $h_{w}$ on the manifolds from $W'{}^{\sigma}$ which contain heteroclinic points with $beh$-type not exceeding $n-1$, and for this manifolds we have constructed trimmed tubular neighborhoods of points on dual manifolds.

Let us extend $h_{g}$ to $h_{w}$ on the manifolds from $W'{}^{\sigma}$ which contain heteroclinic points with $beh$-type, not exceeding $n$. 
Let $W^{\sigma}(z)$ be such a manifold. Choose on $W^{\sigma}(\Omega_{k})$ a fundamental neighborhood  $G^{\sigma}(\Omega_{k})$ such that $\partial G^{\sigma} (\Omega_{k}) \subset W'{}^{u}(f) \cap W '{}^{s}(f)$ and appropriate to it $\widetilde G^{\sigma} (h_{g}( \Omega_{k}))$ on $\widetilde W^{\sigma}(h_{g}( \Omega_{k}))$ such that $h_{g}( \partial G^{\sigma}(\Omega_{k})) = \partial \widetilde G^{\sigma}(h_{g}( G^{\sigma} (\Omega_{k}))$. Put 
$$
U = G^{\sigma}(\Omega_{k}) \setminus \cup_{l} T_{1} U^{\rho, k}(\Omega_{l}) 
$$ 
and, correspondingly,  
$$
\widetilde U =\widetilde G^{\sigma}(h_{g}( \Omega_{k})) \setminus \cup_{l}\widetilde{T_{1}U}^{\rho, k}(h_{g}( \Omega_{l})),
$$
 where $T_{1} U^{\rho, k}(\Omega_{l})$ is a set of trimmed tubular neighborhoods in $W^{\sigma}(\Omega_{k})$ of points from $W^{\rho}(\Omega_{l})$. Using of this sets is correct because those $W^{\rho}(\Omega_{l})$, for which $W^{\rho}(\Omega_{l}) \cap W^{\sigma}(\Omega_{k}) \not= \emptyset $, contain points of $beh$-type not exceeding $n-1$ (otherwise we can show that the maximal  $beh$-type of points on $W^{\sigma}(\Omega_{k})$ must be more than $n$) and we have defined for them the trimmed tubular neighborhoods on the previous steps. According to lemma \ref{lemmaf6}, $\cup_{l} T_{1}U^{\rho, k}(\Omega_{l})$ can be enclosed in a finite number of tubular  neighborhoods $T_{1}U(\gamma_{1})..., T_{1}U (\gamma_{m})$. We extend $h_{w}$ on each of $T_{1}U (\gamma_{i})$ putting in a natural way that $h_{w}(x)=\widetilde \nu_{(\gamma_{i})}(h_{w}(\nu_{(\gamma_{i})}^{-1}(x)))$ ($h_{w}$ has defined on the previous steps). Outside of them, on $U$, we extend $h_{w}$ continuously. 

The following operations (extension $h_{w}$ on the whole $W^{\sigma}(\Omega_{k})$ and construction of  trimmed tubular neighborhoods are the same as in the step 0. Obtained $h_{w}$ is continuous. It can be checked up likewise lemma 2.6 \cite{ps}. \par

\medskip
\noindent
{\bf Part 2.}
 In the previous part we have extended $h_{g}$ to $h_{w}$. To expand $h_{w}$ on the whole manifold, we need trimmed systems of tubular neighborhoods $T_{1}^{\sigma, k}$ and $ \widetilde T_{1}^{\sigma, k}$ for $f$ and $g$ correspondingly. We construct them using already obtained trimmed tubular neighborhoods of heteroclinic points. For each fiber $ t^{\sigma, k}(\xi)$ we obtain its trimmed analog, a fiber $t_1^{\sigma, k}(\xi)$ as a pathcomponent of the set  
$$
t^{\sigma, k}(\xi) \setminus {\bf cl}\Bigl(t^{\rho, k}(h_{w}^{-1}(\widetilde d_{i1}(h_{w}(\xi))) \cup t^{\rho, k}(h_{w}^{-1}(\widetilde d_{i2}(h_{w}(\xi))\Bigr),
$$
containing $\xi$ and obtain $\widetilde t_{1}^{\sigma, k}(\xi)$ as a pathcomponent of 
$$
\widetilde t^{\sigma, k}(\xi) \setminus {\bf cl}\Bigl(\widetilde t^{\rho, k}(h_{w} (d_{i1}(h_{w}^{-1}(\xi))) \cup \widetilde t^{\rho, k}(h_{w} (d_{i2}(h_{w}^{-1}(\xi))\Bigr)
$$
containing $\xi$.
We have used here homeomorphisms $d_{ij}$, defined in lemma \ref{lemmaf1} of this section, specified by the boundary of the tubular neighborhood, which map points of the stable manifold in points of the unstable manifold.

 Let us expand $h_{w}$ from the set ${\bf cl}W'$ on the set $V(\{T_{1}^{\sigma, k}\}) = \cup_{k} V(T_{1}^{\sigma, k})$. According to lemma 1, $V(\{T_{1}^{\sigma, k}\})$ has a structure of direct product. We define $h_w \colon V(\{T_{1}^{\sigma, k}\}) \rightarrow V(\{T_{1}^{\sigma, k}\})$ in the following way: if $x= t_{1}^{\sigma, k}(\xi_{1}) \cap t_{1}^{\sigma, k}(\xi_{2})$ then $h_{v}(x)$ be $t_{1}^{\sigma, k}(\xi_{1}) \cap t_{1}^{\sigma, k}(\xi_{2})$.
The continuity of $h_{v}$ follows from the continuity of $h_{w}$ because of structure of direct product. 

Consider the reminding set $M \setminus V(\{T_{1}^{\sigma, k}\})$. In case when a pathcomponent of the set is contained in a periodic component of the diffeomorphism, the expansion of the $h_{w}$ was described in paper \cite{Bezen}. Consider the case when pathcomponents of the set $M \setminus V(\{T_{1}^{\sigma, k}\})$ are contained in a wandering components of the diffeomorphism. Under the action of $f$ those pathcomponents are divided into a finite number of classes equivalent up to the action of $f$ (it is an easy consequence from the lemma \ref{lemmaf6}).
Every wandering component is homeomorphic to disc \cite{Smith}, and, hence, pathcomponents of the set are homeomorphic to discs.
 Let us choose appropriate representatives $D_{i}$ from each class of equivalence under the action of $f$ and corresponding $\widetilde{D_{i}}$ from each class of equivalence under the action of $g$. We expand $h_{v}$ on them continuously with arbitrary homeomorphisms $k_{D_{i}}$. On the rest of areas we put $h_{f^{n}(D_{i})} = g^{n}\circ k_{D_{i}} \circ f^{-n}$.  

Thus, we have constructed homeomorphism $h: M \rightarrow M $, such that $ h \circ f = g \circ h$.  
It proves the theorem.
\end{proof}

\section{Lattice neighborhoods and local type of heteroclinic points.}
\label{section3}

For the reasoning on induction we assign to every heteroclinic point a numeric value called its local type. 
In this chapter we give a rigorous definition of local type of heteroclinic points, satisfying the lemma \ref{lemma36plus}, and prove other useful properties of the local type.
Note, that construction of the definition of local type of heteroclinic points is specially chosen to satisfy the lemma \ref{lemma36plus}.

\subsection {Definition of lattice structures.}

Let a point $\gamma$ possesses two different lattice tetragons (at that corresponding periodic points may coincide) such that none of the projections of $\gamma$ in one of the tetragons does not belong to the lattice structure corresponding to another tetragon. Let us call them {\em unadjusted tetragons}.

\begin{definition}\label{def2condA}
A set of heteroclinic points which belong to the lattice structure of $p$ satisfies {\em the condition A} if  it does not contain points with unadjusted tetragons. 
\end{definition}

If a set of heteroclinic point satisfies the condition A for $p$ then every point of the set possesses a unique tetragon between this point and $p$. The existence of another tetragon between this point and $p$ leads due to the condition A to homoclinic intersection which is prohibited for Morse-Smale diffeomorphisms.

\begin{definition}\label{def2condB}
A set of heteroclinic points which belong to the lattice structure of $p$ satisfies {\em the condition B} if for every point $\gamma$ of the set if a projection of $\gamma$ on $W^{\sigma}(p)$ belongs to another tetragon with lattice structure then the point $\gamma$ also must belong to the other lattice structure and the projection of $\gamma$ in the other lattice structure along the $\rho$-manifold must belong to the set.
\end{definition}

Let $p$ be a periodic point of $f$.
Consider a set of points that belong to the lattice structure of $p$ and satisfy the conditions A and B.

\begin{lemma}\label{lemma2propAB1}
If the point $\gamma_1$ satisfies the conditions A and B for $p$ then all the points of the tetragon between $\gamma_1$ and $p$ also satisfy the conditions A and B.
\end{lemma}

\begin{proof}
Consider a point $\gamma$ that belongs to the lattice tetragon of $\gamma_1$. 
It is easy to see that $\gamma$ satisfies the condition A  because its other tetragons either contain the point $\gamma_1$, satisfying the condition A, or intersect either $W^{u}(p)$ or $W^{s}(p)$ and $\gamma$ satisfies the condition A directly.

Suppose that the projection of $\gamma$ on $W^{\sigma}(p)$ the point $\pi(\gamma, W^{\sigma}(p))$ belongs to the lattice structure of a periodic point $y$.  If the segment 
$$
[\pi(\gamma, W^{\sigma}(p)), \pi(\pi(\gamma, W^{\sigma}(p)), W^{\rho}(y))] \subset W^{\sigma}(p)
$$
does not contain the point $\pi(\gamma_1, W^{\sigma}(p))$ then we can continue one tetragon into another obtaining that $\gamma$ belongs to the lattice structure of $y$ and the projection $\pi(\gamma, W^{\rho}(y))$ is located inside the tetragon of $\gamma_1$.

Consider the case when the segment 
$$
[\pi(\gamma, W^{\sigma}(p)), \pi(\pi(\gamma, W^{\sigma}(p)), W^{\rho}(y))] \subset W^{\sigma}(p)
$$ 
contains the point $\pi(\gamma_1, W^{\sigma}(p))$. Then $\gamma_1$ also belongs to the lattice structure of $y$ and, according to the condition B, the projection $\pi(\gamma_1, W^{\rho}(y))$ belongs to the set and we also can continue one tetragon into another obtaining that $\gamma$ belongs to the lattice structure of $y$ and the projection $\pi(\gamma, W^{\rho}(y))$ is located inside the tetragon of $\pi(\gamma, W^{\rho}(y))$.
\end{proof}

\begin{lemma}\label{lemma2propAB2}
If the set of points which belong to the lattice structure of $p$ is not empty then the subset of points satisfying the conditions A and B also is not empty.
\end{lemma}

\begin{proof}
It is enough to show the representative of the set. 

Since the set of points which belong to the lattice structure of $p$ is not empty, $W^{u}(p)$ and $W^{s}(p)$ possesses heteroclinic points. At that they must possess points which do not belong to any lattice structure, for example, the points 
      of $beh$-type 1.
Consider a set of points such that its projections on $W^{u}(p)$ and $W^{s}(p)$ are points which do not belong to any lattice structure. They satisfy the condition B by construction. 

Let $\gamma_1$ be the point on $W^{\sigma}(p)$ such that it does not belong to any lattice structure.
Consider the set of heteroclinic points on $\rho$-manifold of $\gamma_1$ which do not satisfy the condition A.  Consider a subset of this set such that every point $\gamma$ of the set belongs to lattice structure of a periodic point $x$ and there exists a heteroclinic point $\zeta_{\gamma} \in W^{\rho}(x)$ which does not belong to lattice structure of $p$ and is a projection of $\gamma$ on $W^{\rho}(x)$. Since $W^{\rho}(x)$ cannot be a boundary for itself, all points of the set are contained in a segment of  $W^{\rho}(x)$. Furthermore, the set is closed as the intersection of lattice tetragons. If the boundary point of the set is not contained on $W^{\sigma}(p)$ we can choose a new point $\gamma_2$ such that both its projections  have the $beh$-type 1 and its tetragon does not contain this set. In case when $W^{\sigma}(p)$ contains the boundary point we have a contradiction, since in this case all points of the set belong to the lattice structure of $p$.

As the amount of periodic points is finite, we can choose point $\gamma^*$ which satisfy the conditions A and B. 
\end{proof}

\begin{lemma}\label{lemma2propABClose}
A set of points which belong to the lattice structure of a periodic point $p$ and satisfy the conditions A and B is a closed set provided we added to the set the point $p$ and heteroclinic points on $W^{u}(p)$ and $W^{s}(p)$.
\end{lemma}

\begin{proof}
It is easy to see that the limit point of a sequence of points that belong to the lattice structure of $p$ also belongs to the lattice structure of $p$. Let us show that the limit point satisfies the condition B. 

Suppose on the contrary that there exists a limit point $\gamma \in W^{\sigma}(x) \cap W^{\rho}(y)$ such that it does not belong to the lattice structure of a periodic point $p_1$ but its projection $\gamma_2$ on $W^{\sigma}(p)$ does. Since $\gamma_2$ is the limit point of  projections of points of the sequence, we can choose a subsequence whose points belong to the lattice structure of $y$ and, hence, its projections belong to the lattice structure of $p_1$. Since they satisfy the condition B, they also belong to the lattice structure of $p_1$ and $\gamma$ as its limit point also belongs to the lattice structure of $p_1$. 

Let us show that the limit point satisfies the condition A. 
Suppose on the contrary that there exists a limit point $\gamma \in W^{\sigma}(x) \cap W^{\rho}(y)$ such that it belongs to the lattice structure of a periodic point $p_1$ and the corresponding tetragons are unadjusted. Since $W^{\rho}(y)$ has a neighborhood where all heteroclinic points belong to its lattice structure, there is a subsequence of points converging to $\gamma$ and belong to the lattice structure of $y$. Continuing the tetragons, we have that the points of the subsequence belong to the lattice structure of $p_1$. At that they satisfy the condition A. We have a contradiction.
\end{proof}

\subsection {Definition of lattice neighborhoods.} \par

For a heteroclinic point $\gamma \in W^{\sigma}(z) \cap W^{\rho}(p)$ consider the set of points which belong to the lattice structure of $p$, satisfy the conditions A and B and $\gamma$ is their projection on $W^{\rho}(p)$ in the lattice structure of $p$. It is easy to see from lemmas \ref{lemma2propAB1}, \ref{lemma2propAB2}, \ref{lemma2propABClose} that there exist heteroclinic points $\gamma^{*}_1, \gamma^{*}_2 \in W^{\sigma}(z)$ such that the set defined above is exactly the set of heteroclinic points which belong to the connected segment $[\gamma^{*}_1, \gamma^{*}_2] \subset W^{\sigma}(z)$ assuming that $\gamma$ also belongs to the set. 
Let us call it {\em the neighborhood of $\gamma$ with conditions A and B}.

\begin{lemma}\label{lemma2nbcndAB}
if neighborhoods with conditions A and B of points  $\gamma_1 \in W^{\sigma}(z) \cap W^{\rho}(p)$ and $\gamma_2 \in W^{\sigma}(z) \cap W^{\rho}(y)$ have non-empty intersection then one of the points is contained in the neighborhood of another point.
\end{lemma}

\begin{proof}
The case when the points are contained in the neighborhoods of each other is impossible here since it leads to homoclinic intersection.

Let, for example, $\gamma_1$ is not contained in the neighborhood with conditions A and B of the point $\gamma_2$. Let $\gamma$ be the point of the intersection. $\gamma$ belongs to the neighborhood of $\gamma_1$, the corresponding tetragons are adjusted, so $\pi(\gamma, W^{\sigma}(p))$ belongs to the lattice structure of $y$. Then, according to the condition B for $p$,  $\gamma_2$ which is $\pi(\gamma, W^{\rho}(y))$, belongs to the neighborhood of $\gamma_1$.
\end{proof}

\medskip
\smallskip
Consider the neighborhood with conditions A and B of point  $\gamma_1 \in W^{\sigma}(z) \cap W^{\rho}(p)$. 
$W^{\sigma}(p)$ contains a countable set of points which do not belong to any lattice structure. Being the projections of points from the neighborhood with conditions A and B of point  $\gamma_1$, they separate neighborhoods of points which are contained in the neighborhood with conditions A and B of point  $\gamma_1$, moreover, if a point $\gamma_2$ is contained in the neighborhood of $\gamma_1$ and there exists a point $\gamma^*$ such that $\pi(\gamma^*, W^{\sigma}(p))$ does not belong to the lattice structure of any point and $\gamma_2$ belongs to the segment $[\gamma_1,\gamma^*] \subset W^{\sigma}(z)$ then the neighborhood with conditions A and B of $\gamma_2$ is entirely contained in the neighborhood  with conditions A and B of $\gamma_1$. This property may not take place only for neighborhoods which intersect the neighborhood of $\gamma_1$ after the last point like $\gamma^*$, near the boundary of the neighborhood.

Consider the set of heteroclinic points on $W^{u}(p)$ and $W^{s}(p)$ with the property that they do not belong to neighborhood with condition A and B of any periodic point. 
Let us call it {\em the sets F}.

For a point $\gamma_0 \in W^{\sigma}(p)$ of the set F let us define {\em the wide neighborhood of $\gamma_0$ in $W^{\sigma}(p)$} as a set of heteroclinic points with the following property: if a point $\gamma$ belongs to the wide neighborhood of $\gamma_0$ then there exists a sequence of points $\gamma =\gamma_{n}, \gamma_{n-1}, \ldots, \gamma_0$ such that $\gamma_i$ belongs to the neighborhood with conditions A and B of  $\gamma_{i-1}$. 

\begin{lemma}\label{lemma2prWN}
Wide neighborhoods of different points $\gamma_0$ and $\gamma_0^*$ do not intersect one another.
\end{lemma}

\begin{proof}
Suppose on the contrary that they have a common point $\gamma$. Then, according to the definition, there exists a sequence of points $\gamma =\gamma_{n}, \gamma_{n-1}, \ldots, \gamma_0$ such that $\gamma_i$ belongs to the neighborhood with conditions A and B of  $\gamma_{i-1}$ and a sequence of points $\gamma =\gamma^*_{m}, \gamma_{m-1}, \ldots, \gamma^*_0$ such that $\gamma^*_i$ belongs to the neighborhood with conditions A and B of  $\gamma^*_{i-1}$
Since $\gamma$ is a common point of neighborhoods of  $\gamma_{n-1}$ and $\gamma^*_{m-1}$ then one of the points must belong (see lemma \ref{lemma2nbcndAB}) to the neighborhood of another. Let, for example, $\gamma_{n-1}$ is contained in the neighborhood of $\gamma^*_{m-1}$. Then the points $\gamma_{n-2}$ and $\gamma^*_{m-1}$ have the common point $\gamma_{n-1}$ and so on. Finally, we obtain that $\gamma_{0}$ and $\gamma^*_{0}$ must have the common point. We obtain the contradiction.
\end{proof}

\medskip
Introduce the relation of equivalence on the set of points belonging to the lattice structure of $p$ and satisfying the conditions A and B. Two heteroclinic points are equivalent (belong to the same {\em wide tetragon}) if and only if their projections on $W^{u}(p)$ and $W^{s}(p)$ belong to the same wide neighborhoods.

\begin{definition}
\label{definition201}
A lattice neighborhood of  $p$ is the point $p$ and a set of heteroclinic points which either belong to $W^{u}(p)$ and $W^{s}(p)$ or belong to the subset of the set of points belonging to the lattice structure of $p$ and satisfying the conditions A and B with the additional property that a point is contained in the subset if and only if its wide tetragon is entirely contained in the set.
\end{definition}


It is easy to see that the lattice neighborhood preserves the properties of the set of points with conditions A and b:
 if the lattice neighborhood of a periodic point $p$ possesses a heteroclinic point $\gamma$ then it possesses all heteroclinic points that are situated in tetragon with lattice structure between $p$ and $\gamma$, lattice neighborhood is a closed set and so on.

Though a lattice neighborhood is a countable set of heteroclinic points we call it  ``neighborhood'' since it can be enclosed in an invariant under the action of $f^{n}$  neighborhood of $p$ that contains $W^{u}(p)$ and $W^{s}(p)$ to look like the tubular neighborhood of $p$.

\medskip
Let $\gamma$ be a heteroclinic point of the intersection of $W^{u}(x)$ and $W^{s}(y)$. 
\begin{definition}
{\em The lattice neighborhood of a heteroclinic point $\gamma$} is a set of heteroclinic points that belong to the intersection of $W^{u}(x)$ ($W^{s}(y)$) with the lattice neighborhood of $y$ ($x$) such that $\gamma$ is their projection on $W^{s}(y)$ ($W^{u}(x)$).
\end{definition}

\begin{definition}
The union  of lattice neighborhoods of all points of a heteroclinic trajectory is called {\em a lattice neighborhood of the heteroclinic trajectory}. 
\end{definition}

\begin{remark}\label{remark28}
if two lattice neighborhoods of heteroclinic points on the same manifold have non-empty intersection then one of them contains the other. 
\end{remark}

This property follows from the procedure of construction of lattice neighborhoods.

\medskip

\subsection {Definition of local type of points.}

\begin{definition}
{\em
a heteroclinic point  on $W^{\sigma}(x)$ is called {\em a heteroclinic point of local type 1 on $W^{\sigma}(x)$} if this point is not contained in the lattice neighborhood on $W^{\sigma}(x)$ of any other heteroclinic point;
 it is called {\em a point of local type 2 on $W^{\sigma}(x)$} if on $W^{\sigma}(x)$ this point is contained only in lattice neighborhoods of heteroclinic points of the first local type; 
it is called {\em a point of local type $n$ on $W^{\sigma}(x)$} if on $W^{\sigma}(x)$ this point is contained only in lattice neighborhoods of heteroclinic points up to $n-1$th local type inclusively. \par

The local type of a  point $\gamma$ on $W^{\sigma}(x)$  is denoted by $loctype (\gamma, W^{\sigma}(x))$.
}
\end{definition}

\begin{remark}
Since the lattice neighborhoods are invariant under the action of $f$, all points that belongs to the same trajectory have the same local type. 
\end{remark}

\begin{lemma}\label{lemma34}
The manifold contains a finite number of trajectories of local type 1. 
\end{lemma}

\begin{proof}
 It is enough to prove, that there exists a finite number of points with local type 1 on each fundamental neighborhood. 
Let us assume by contradiction that there exist infinitely many points of local type 1 on a fundamental neighborhood of a manifold. Then the set of this points has at least one boundary point, which, according to \ref{intr13}-1 is heteroclinic. We obtained a contradiction because all points which lay in the lattice neighborhood of this boundary point cannot have local type 1 according to the definition of local type. \par
\end{proof}

\begin{lemma}\label{lemma35}
Let $\gamma_1$ be a heteroclinic point from $W^{\sigma}(x) \cap W^{\rho}(z)$ and $\gamma \in W^{\sigma}(x) \cap W^{\rho}(y)$ be a heteroclinic point from the lattice neighborhood of the point $\gamma_1$. 
Then 
$$
loctype (\gamma, W^{\sigma}(x)) = loctype (\gamma_1, W^{\sigma}(x)) + loctype (\pi (\gamma), W^{\sigma}(z)),
$$
 where $\pi (\gamma)$ is a projection of the point $\gamma$ along the manifold $W^{\rho}(y)$ on the manifold $W^{\sigma}(z)$. 
\end{lemma}

\begin{proof}
The lemma is illustrated by figure \ref{fig1}.
According to the definition, $loctype (\gamma, W^{\sigma}(x))$ shows the number of points, whose lattice neighborhoods contains the point $\gamma$. Let us divide them into two classes: points which lay in lattice neighborhood of $\gamma_1$ and points which do not lay in but contain in their lattice neighborhoods the point $\gamma_1$. The number of points in the second class is given by $loctype (\gamma_1, W^{\sigma}(x))$. 

Points of the first class can be projected on $W^{\sigma}(z)$ using the lattice structure.
Since their lattice neighborhoods on  $W^{\sigma}(x)$ contain $\gamma$ all the more their projections contain in their lattice neighborhoods on $W^{\sigma}(z)$ the point $\pi (\gamma)$. 
On the contrary, If $\pi (\gamma)$ lays in a lattice neighborhood of a point $\pi (\zeta)$ then the definition of lattice neighborhoods demands $\rho$-manifold of this point to intersect $W^{\sigma}(x)$ in a point $\zeta$ so that $\gamma$ belongs to lattice neighborhood of $\zeta$ and $\zeta$ belongs to lattice neighborhood of $\gamma_1$. Therefore the number of points in the first class is given by $loctype (\pi (\gamma), W^{\sigma}(z))$. It proves the lemma. 
\end{proof}

\begin{lemma}\label{lemma36}
Local type of a heteroclinic point on its stable manifold is equal to its local type on the unstable manifold. \par
\end{lemma}

\begin{proof}
Consider heteroclinic point $\gamma \in W^{\sigma}(x) \cap W^{\rho}(y)$. There are two possible cases: first case when one of local types of the point is equal to 1 and second case when both local types of the point is greater than 1.
Consider the first case. Let the point have the local type 1 on its $\rho$-manifold, $loctype (\gamma, W^{\rho} (y)) = 1$. Assume by contradiction that $loctype (\gamma, W^{\sigma}(x)) = k> 1$. Then, according to the definition of local type, there exists a point $\gamma_{1} \in W^{\sigma}(x) \cap W^{\rho}(z)$ of local type $k-1$, which contains in its lattice neighborhood the point $\gamma$. 
Since $\gamma$ belongs to the lattice neighborhood of the periodic point $z$, then, according to definition, $W^{\rho}(y)$ intersects $W^{\sigma}(z)$ in a point $\gamma_{2}$ (see fig. 1). 
Consider obtained tetragon with lattice structure, formed by points $z$, $\gamma_{1}$, $\gamma$ and  $\gamma_{2}$ from another point of view as a tetragon of the point $\gamma$ on $W^{\rho}(y)$. We have, according to definition, that the point $\gamma$ is contained in the lattice neighborhood of the point $\gamma_{2}$ and, hence, $loctype (\gamma, W^{\rho}(y))>1$. We have obtained a contradiction. \par

In the second case we prove the lemma on induction. The first case gives us the base of induction.
Let for heteroclinic points for which both local types is less $k$, the lemma holds true. Let us prove it for a point $\gamma$, for which both its local types are less $k + 1$. We make use of  lemma \ref{lemma35} because both projections of $\gamma$ exist. Let $\gamma$, $\gamma_{1}$, $\gamma_{2}$ be such as in fig. 1. We have the following relations: 
$$k + 1>Loctype (\gamma, W^{\sigma}(x))  = loctype (\gamma_{1}, W^{\sigma}(x)) + loctype (\gamma_{2}, W^{\sigma}(z)$$
$$k + 1> Loctype (\gamma, W^{\rho}(y)) = loctype (\gamma_{2}, W^{\rho}(y)) + loctype (\gamma_{1}, W^{\rho}(z)) $$
Since the local type of a heteroclinic point is a positive number, both local types of points $\gamma_{1}$ and $\gamma_{2}$ are less than $k$, and, according to the induction supposition, 
$$loctype (\gamma_{1}, W^{\sigma}(x)) = loctype (\gamma_{1}, W^{\rho}(z))$$
$$loctype (\gamma_{2}, W^{\sigma}(z)) = loctype (\gamma_{2}, W^{\rho}(y))$$ 
and we have $loctype (\gamma, W^{\sigma}(x)  = loctype (\gamma, W^{\rho}(y))$. \par
\end{proof}

   Thus, we define \textit{local type of a heteroclinic point} as local type of a point on one of its manifolds. Let us denote local type of a heteroclinic point of $\gamma{}$ as $loctype (\gamma) $. \par


\begin{lemma}\label{lemma36plus}
the local type of a heteroclinic point is equal to sum of local types of its projections in any lattice neighborhood to which the point belongs.
\end{lemma}

\begin{proof}
It is a consequence of lemmas \ref{lemma35} and \ref{lemma36}.
\end{proof}

\begin{lemma}\label{lemma37}
Consider  the lattice neighborhood  of a point of  local type $k-1$. If it contains heteroclinic points of local type $k$ then their set is infinite and it has no other boundary points except this point of local type $k-1$. 
\end{lemma}

\begin{proof}
Let $\gamma_{1} \in W^{\sigma}(x) \cap W^{\rho}(z)$ be a heteroclinic point of the local type $k-1$ and $\gamma \in W^{\sigma}(x) \cap W^{\rho}(y)$ be a heteroclinic point of the local type $k$, laying in the lattice neighborhood of $\gamma_{1}$ on $W^{\sigma}(x)$. According to the definition of the lattice structure, $W^{\rho}(y)$ intersects $W^{\sigma}(z)$ in a point $\gamma_{2}$, and there appears a tetragon with lattice structure formed by points $z$, $\gamma_{1}$, $\gamma{}$, $\gamma_{2}$.
In this tetragon $W^{\rho}(y)$ intersects $W^{\sigma}(z)$ infinitely many times in points of the trajectory $\gamma_{2}$, and, hence, infinitely many times intersects $W^{\sigma}(x)$. Using the lemma \ref{lemma35}, one readily sees that points of the trajectory $\gamma_{2}$ have the local type 1, and the appropriate points on $W^{\sigma}(x)$ have the local type k-1. \par

If the set of heteroclinic points of local type $k$ has anther boundary point then due to the lattice structure the projection of this boundary point is a boundary point for the set of projections of local type 1 of points of local type $k$. But it is impossible because this set has a unique boundary point $z$ which is the projection of $\gamma_1$ on $W^{\sigma}(z)$.
\end{proof}

\section {Algorithm of construction of the extended coding set. }
\label{Sec5AlgCS}

This section is devoted to the proof of the theorem \ref{TheoremCS2} (see section \ref{Sec4defCS}).

\begin{vtheorem}{theorem}{\ref{TheoremCS2}}
The basic coding set and the graph $G(f)$ (both finite) are uniquely determine the extended coding set. The obtained extended coding set is also finite.
\end{vtheorem}

\begin{proof}
We construct the extended coding set from the basic coding set using induction on local type of points\footnote{We also use the information on orientation and incidence in periodic points of the diffeomorphism, stored in a finite graph $G(f)$ which is defined in section \ref{Sec7graph}.}.

To describe the trajectories of local type 1, we simply add their simple formulas to the set. Since the step of induction in general is too complicated, it is quite appropriate here to describe first as an example the case of points of local type 2 before considering a general case.

\medskip
{\bf Proof. The extended coding set for points of local type 2.}

On the previous step we obtain a finite set of  formulas in the extended set for points of local type does not exceeding 1.
 
 Now we modify every formula of extended coding set using formulas of basic coding set and  graph  to obtain new additional formulas that describes trajectories of local type 2. After obtaining the formulas to proceed with inductive construction we construct an auxiliary technical set which is used in the process of further construction . We have not constructed it on the previous step of induction because for that case it is directly contained in the basic coding set.

\begin{figure}[htbp]
\includegraphics{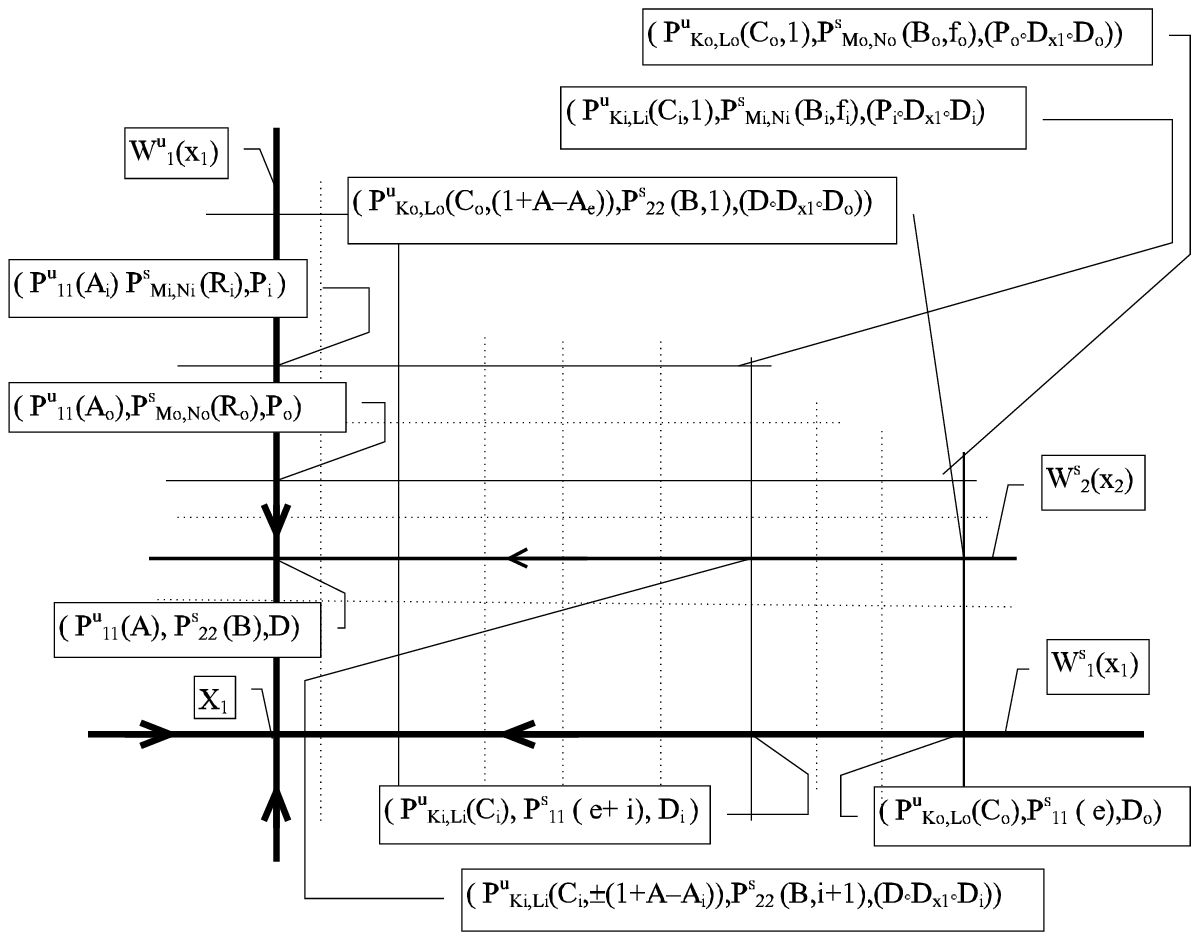}
\caption{Illustration to examples of this section.}\label{fig2}
\end{figure}

For example, (see fig. \ref{fig2}) let us describe all trajectories of local type 2, laying in the back set of the lattice neighborhood of a point $(P^{u}_{11} (A), P^{s}_{22}(B), D)$, where extreme point is a point  $(P^{u}_{K_0, L_0}(C_{o}, A-A_{o}), P^{s}_{22}(B, 1), D_{o}* D * D_{x1}\protect\footnote {Here $D$ is the orientation of the point ($P^{u}_{11}(A), P^{s}_{22}(B), D)$, $D_{o}$ is the orientation of the point $( P^{u}_{K_0, L_0}(C_{o}), P^{s}_{11}(e), D_{o}), D_{x1}$ is the   orientation of a frame $W^{u}_{1}(x_{1})$ --- $x_{1}$ --- $W^{s}_{1}(x_{1})$. First two values is contained in the basic coding set, the third one is contained in the graph. The orientation of the extreme point is equal to the product of these three values. Also the second component of the point number on $W^{u}_{K_0}(x_{ L_0})$ is substituted by the obvious expression. It is done to decrease the number of designations.})$.

As mentioned above,  the order of  disposition of points of local type 2 in a lattice neighborhood repeats the order of  disposition of their projections on $W^{s}_{1}(x_{1})$ which are points of local type 1.
Comparing the orientation of $D$ with a mutual location of manifolds $W^{u}_{1}(x_{1})$, $W^{s}_{1}(x_{1})$, $W^{s}_{2}(x_{1})$, stored in the graph, we determine, for example, that a projection of the back set on $W^{s}(x_{1})$ lays, in fact, on $W^{s}_{1}(x_{1})$, and directions of the numeration of the back set and its projection coincide. Thus, points of the projection should look like $(P^{u}_{K_0, L_0}(C_{o})\ldots\quad $. 
After searching the basic coding set we obtain a simple formula of this point. For example, let a projection of the extreme point be a point $(P^{u}_{K_0, L_0} (C_{o}), P^{s}_{11}(E), D_{o})$. 

A fundamental neighborhood $[ P^{s}_{11}(E), P^{s}_{11}(E + N_{11}^{s}) )^{s} \subset W^{s}_{1}(x_{1})$ contains $N_{11}^{s}$ of points of local type 1. Let they be the points 
$$
(P^{u}_{K_i, L_i}(C_{i}), P^{s}_{11}(E + i), D_{i})_{i = 0... N_{11}^{s} - 1}
$$
 For each of these points the extreme coding set now contains simple formulas of its trajectories. Furthermore, the  basic coding set contains simple formulas of extreme trajectories of local type 2 for lattice neighborhoods of every trajectory of local type 1. 
Consider the extreme points of lattice neighborhoods of points $(P^{u}_{K_i, L_i}(C_{i}), P^{s}_{11}(E + i), D_{i})$, laying on the same side of  $W^{s}(x_{1})$ where is the manifold $W^{s}_{1}(x_{1})$. Let they be the points 
$$
(P^{u}_{K_i, L_i}(C_{i}, \pm 1\protect\footnote {a sign $+$ or $-$ depends in formulas whether this set is the back one or the front one.}), P^{s}_{M_i, N_i} (B_{i}, f_{i}), P_{i }* D_{x1}* D_{i})_{i = 0... N_{11}^{s} - 1}
$$
 The first component of the prefix of the point's number on its stable manifold shows us the projection of this extreme point on $W^{u}_{1}(x_{1})$ which is a point of local type 1 according to the properties of local type. Let the projections of these extreme points on $W^{u}_{1}(x_{1})$ be points 
$$
(P^{u}_{11} (A_{i}) P^{s}_{M_i, N_i}(R_{i}), P_{i})_{i = 0... N_{11}^{s} - 1}
$$
 
Now we have enough information to start construction of extended coding set formulas that describe points of this back set. Using obvious calculations, we obtain that first $N_{11}^{s}$ of points of local type 2 in our neighborhood can be noted as 
$$
(P^{u}_{K_i, L_i} (C_{i}, \pm (1 + A - A_{i})), P^{s}_{22} (B, i + 1), D * D_{x1}* D_{i})_{i = 0... N_{11}^{s} - 1}
$$
 Since the set is invariant, the whole trajectory of a extreme point consist of extreme points. Hence, on the next fundamental neighborhoods $[ P^{s}_{11} (E + k * N_{11}^{s}), P^{s}_{11} (E + (k + 1)* N_{11}^{s}))^{s} \subset W^{s}_{1}(x_{1})$, $k \geq 0$, the extreme points are 
$$
(P^{u}_{K_i, L_i} (k * N^{s}_{K_i, L_i} + C_{i}, \pm 1), P^{s}_{M_i, N_i} (k * N^{s}_{M_i, N_i} + B_{i}, f_{i}), P_{i}* D_{x1}* D_{i})_{i = 0... N_{11}^{s} - 1,\ k \geq 0}
$$ 
Similarly, all points of local type 2, laying in the back set of the lattice neighborhood of the point $(P^{u}_{11} (A), P^{s}_{22} (B), D)$, are described with formulas 
$$
\left(
\begin{array}{ccc}
P^{u}_{K_i, L_i} (C_{i} + k * N^{u}_{K_i, L_i}, -D_{i}* (1 + |A - A_{i} - k * N^{u}_{11}|)), \cr P^{s}_{22}(B, i + k * N^{s}_{22}+1), \cr
D * D_{x1}* D_{i}
\end{array}
\right)_{ k \geq 0,\ i = 0... N_{11}^{s} - 1}
$$

Obtained formulas do not describes the trajectories yet but describes their representatives. Since trajectories are obtained from points under the action of the diffeomorphism and lattice neighborhoods are invariant, to describe trajectories we simply add a displacement under the action of the diffeomorphism in the first component of point's prefix. We have that these trajectories are described by formulas 
$$
\left(
\begin{array}{ccc}
P^{u}_{K_i, L_i}
\left(
C_{i} + (t + k)* N^{u}_{K_i, L_i},
 -D_{i}* (1 + |A - A _{i} - k * N^{u}_{11}|
\right), \\
P^{s}_{22 } (B + t * N^{s}_{22}, i + k * N^{s}_{11} + 1),\\ 
D * D_{x1}* D_{i}
\end{array}
\right)_{k \geq 0,\ t \in \zz, \atop i = 0... N_{11}^{s}-1}
$$

Repeating this operations for all trajectories of local type 1, we receive in the end the extended coding set for points of local type 2 which is a finite union of the formulas.

\medskip
{\bf Proof. The auxiliary coding set for points of local type 2.}

On the next steps of induction we will need a technical auxiliary set of formulas associated with obtained formulas for points of local type 2, which have been constructed here.
This set of formulas contains two subset of formulas.

The first subset contains formulas associated with points of local type 2, determining extreme points of all local types from 3rd up to maximal  one.
When we substitute in a formula concrete values of its parameters we obtain a heteroclinic point. In our case, it is a point of local type 2. The fact that a formula is associated means that the associated formula have the same domain of definition as the formula to which it is associated and when it substituted with the same value of parameters it gives an extreme point  of the lattice neighborhood of a point of local type 1 such that
it lays in the same lattice neighborhood of a point of local type 1 as a corresponding point of local type 2 and on the same side.

The second subset contains associated formulas determining points of local type 1 on dual manifolds to  manifolds which contain in their lattice neighborhoods these extreme points.

The idea of their construction is very simple:
Extreme points are described in the basic coding set. In formulas the first components of vectors contain information about the point of local type 1 to whose lattice neighborhood extreme points belong. Similarly, the information about the points of local type 1 is contained in the first components of vectors in formulas describing points of local type 2. Substituting their functions that calculate the first component of vector instead of the trajectory parameter we obtain the required auxiliary set.

To give an example of the construction, let us construct formulas, describing extreme points of local type $n$, associated to the formulas of the previous example where we described points of local type 2. 
This example is illustrated by fig. \ref{fig2}.

Let us take $N_{11}^{s}$ of formulas from the basic coding set describing extreme points of local type $n$ for trajectories 
$$
(P^{u}_{K_i, L_i} (C_{i}+t* N^{u}_{K_i, L_i}), P^{s}_{11} (E+i +t* N^{s}_{11}), D_{i})_{  t \in \zz,\ i = 0... N_{11}^{s} - 1}
$$
Trajectories of this extreme points look like 
$$
(P^{u}_{K_i, L_i}(t * N^{s}_{K_i, L_i} + C_{i}, \underbrace{\pm 1, \ldots , \pm 1}_{ n-1\rm\; times}), ***, ***
\protect\footnote {here $***$ substitutes inessential concrete aspect of the part of the formula.}
)_{ t\in \zz,\ i = 0... N_{11}^{s} - 1}
$$
The formulas of trajectories of local type 2 are:
$$
\left(
P^{u}_{K_i, L_i}
\left(
C_{i} + (t + k)* N^{u}_{K_i, L_i},\ ***
\right), 
***,\ ***
\right)_{ k \geq 0,\ t \in \zz,\ i = 0... N_{11}^{s}-1}
$$
Substituting $t+k$, where $k \geq 0,\ t \in \zz$, instead of $t$ in the formulas of extreme trajectories, we obtain required associated formulas
$$
(P^{u}_{K_i, L_i}((t+k) * N^{s}_{K_i, L_i} + C_{i}, \underbrace{\pm 1, \ldots , \pm 1}_{ n-1\rm\; times}), ***, ***)_{ t\in \zz,\ i = 0... N_{11}^{s} - 1}
$$.

From another point of view, we applied to the formulas of extreme points the same operations as we performed 
constructing the formulas of trajectories of local type 2, namely, 
1)limitation of the trajectory parameter (by substituting $t = k$, where $k \geq 0$; 
2)adding a displacement under the action of a diffeomorphism, depending of a new trajectory parameter $t\in\zz$. 
We receive that constructed  formulas for points of local type 2 and extreme points depend on the same parameters. 

In the same way we construct associated formulas of second kind, describing points of local type 1.
This points contain in their lattice neighborhoods the extreme points of local type $n$ and are their projections through the manifolds dual to manifolds which contain associated trajectories of local type 2.

For example, let the points of local type 2 and the extreme points be the same as in the previous example. Then the points of local type 1 are determined uniquely by the prefix and the first component of vector of the $s$-part of the associated formula for extreme points of local type $n$ and we can find them searching the basic coding set for their simple formulas. We obtain the required auxiliary formulas substituting the functions that calculate the first component of vector in the extreme points instead of the trajectory parameter.

Thus, on the second step of induction we have constructed the extended coding set containing the explicit formulas of trajectories of local types up to the second one, and the auxiliary coding set required to simplify further construction. All the sets are finite.

\medskip
{\bf Proof. The general step of induction.}

Suppose the extended coding set containing finite set of formulas for trajectories up to $(n -1)th$ local type inclusively, where formulas for trajectories of local type $k$ depend from trajectory parameter and $k-1$ of local parameters, and auxiliary coding set of associated formulas is already constructed. 
Let us add to the extended coding set formulas for trajectories of local type $n$.

The idea is the same as in the previous step.
Consider a lattice neighborhood on $W^{\sigma}(z)$ of a trajectory of local type 1 belonging to $W^{\sigma}(z) \cap W^{\rho}(x)$. To describe trajectories of local type  $n$ laying in the lattice neighborhood we first consider their projections on $W^{\sigma}(x)$.

The order of  location of points of local type $n$ in a lattice neighborhood repeats the order of  location of their projections on $W^{\sigma}(x)$ which are points of local type $n-1$ according to lemma \ref{lemma35}. To construct the $\sigma$-part of the required formula we paste the prefix and the first component of the vector from the $\sigma$-part of  the simple formula of the heteroclinic trajectory of local type 1 with vectors from every  $\sigma$-part of the formulas describing projections and substitute a limited parameter for the trajectory parameter. (We require the limited parameter because formulas with the trajectory parameter also describe points which are not projections). 

Then we calculate the $\rho$-part of the required formulas. Their prefixes and $n-1$ of first components of vectors are exactly the same as in the associated formulas that describe extreme trajectories of local type $n$. The modulo of value of the last component of vector is equal to the number of heteroclinic points of local type 1 situated between our point of local type 1 and projection of local type 1 of the associated extreme point of local type $n$, stored in the auxiliary coding set. The orientation is calculated as product of the orientation in the formula for points of local type $n-1$, the orientation of point of local type 1 and the orientation in the periodic point $x$.

Let us show an example of the construction. Let us describe all trajectories of local type $n$, laying in the back set of the lattice neighborhood of the point $(P^{u}_{11} (A), P^{s}_{22}(B), D)$. 
This example is also partially illustrated by fig.~\ref{fig2}.
Projections of the trajectories of local type $n$ on $W^{s}_{1}(x_{1})$  have the local type $n-1$. 
According to the induction supposition, the projections on $W^{s}_{1}(x_{1}) $ of local type $n-1$ are described by a finite set of formulas depending from the trajectory parameter and $n-2$ of local parameters. Let us restrict domain of definition of parameters so that these formulas will describe only projections of the trajectories of local type $n$, imposing a finite number of restrictions. We do it in the following way: let extreme point of our set is projected in the lattice neighborhood on $W^{s}_{1}(x_{1})$ of a point of local type 1, corresponding to a value $t = t_{0}$ of the trajectory parameter. Then for every  value $t>t_{0}$ of the trajectory parameter all obtained points are projections, and for every  value $t < t_{0}$ are not projections  and we add a restriction $t\in [t_{0}+ 1, +\infty) \cap \zz$ (at that the trajectory parameter turns in local one). When $t =t_{0}$ in the same way we construct restriction for the next local parameter $k_{2}$ and so on. These new restrictions we also apply to the associated formulas. 

Now we transform each formula for points of local type $n-1$ into the formula for points of local type $n$ in the lattice neighborhood of a point $(P^{u}_{11} (A), P^{s}_{22} (B), D)$. We obtain $s$-part of the required formula, pasting $\{P^{s}_{22}(B, \}$--- (the prefix and the first component of the $s$-part of this point) with the vectors of $s$-parts of formulas. The $u$-part of the required formulas is constructed in the following way: the prefix and $n-1$ first components we get from associated formulas for extreme points of local type $n$.
The $n$th component is computed as $(-D^*)* (1 + |A-A(k_{1}(t), k_{2}..., k_{n-1}) |)$, where A is taken from the first component of the $u$-part of  $(P^{u}_{11} (A)\dots$, function of orientation $D^* = D^*(k_{1}(t), k_{2}..., k_{n-1})$ is taken from associated formula for extreme points of local type $n$ and $(k_{1}(t), k_{2}..., k_{n-1})$ is taken from the first component of the $u$ - part of associated formula for projections of local type 1 of extreme points of  local type $n$ on $W^{u}_{1}(x_{1})$. 
The orientation of the required formula is calculated as $ D * D^{*} * D_{x_1}$. 
Passing from points to trajectories, we introduce a trajectory parameter $t$ and add a displacement under the action of the diffeomorphism, corresponding to this parameter,  in the first components of  the $u$- and $s$- parts of the formulas. The set of obtained formulas is finite. 

If our lattice neighborhood of  $(P^{u}_{11} (A), P^{s}_{22} (B), D)$ contains points of local type, greater than $n$, we also need to construct  associated formulas to the already obtained ones for extreme points of greater types and their projections of local type 1, which will be used on the next step of induction. This is a scheme of construction of associated formulas for extreme points of local type $k$: 

According to the induction supposition,  such associated formulas exist for formulas of points of local type $n-1$. Let us apply to them the same transformations as applied to formulas of trajectories of local type $n-1$ to obtain formulas of trajectories of local type $n$. First we introduce the restricted domain of definition of parameters , then add the same displacement under the action of the diffeomorphism.
We receive the $n$-parametrical formulas, describing extreme points of local type $k$, associated to the obtained formulas of trajectories of local type $n$. 

The prefix and the first component of vector of the corresponding part of the associated formula for extreme points of local type $k$ determines their projections of local type 1. Their associated formulas are obtained by substituting the functions that calculate the first component of vectors of the extreme points in simple formulas of projections of local type 1 instead of their trajectory parameter.
 
We have constructed the extended coding set containing a finite number of formulas and describing trajectories of local types up to $n$th one and the auxiliary coding set of associated formulas.

When $n$ is equal to highest possible local type, we obtain the required finite extended coding set.
\end{proof}

\end {document}